# An efficient Hessian based algorithm for solving large-scale sparse group Lasso problems

Yangjing Zhang[*]   Ning Zhang[†]   Defeng Sun [‡]   Kim-Chuan Toh [§]

December 15, 2017

**Abstract.** The sparse group Lasso is a widely used statistical model which encourages the sparsity both on a group and within the group level. In this paper, we develop an efficient augmented Lagrangian method for large-scale non-overlapping sparse group Lasso problems with each subproblem being solved by a superlinearly convergent inexact semismooth Newton method. Theoretically, we prove that, if the penalty parameter is chosen sufficiently large, the augmented Lagrangian method converges globally at an arbitrarily fast linear rate for the primal iterative sequence, the dual infeasibility, and the duality gap of the primal and dual objective functions. Computationally, we derive explicitly the generalized Jacobian of the proximal mapping associated with the sparse group Lasso regularizer and exploit fully the underlying second order sparsity through the semismooth Newton method. The efficiency and robustness of our proposed algorithm are demonstrated by numerical experiments on both the synthetic and real data sets.

[*]Department of Mathematics, National University of Singapore, 10 Lower Kent Ridge Road, Singapore 119076 (`zhangyangjing@u.nus.edu`).

[†]Department of Mathematics, National University of Singapore, 10 Lower Kent Ridge Road, Singapore 119076. Future Resilient Systems, Singapore-ETH center. The Singapore-ETH Centre (SEC) was established as a collaboration between ETH Zurich and National Research Foundation (NRF) Singapore (FI 370074011) under the auspices of the NRF's Campus for Research Excellence and Technological Enterprise (CREATE) programme (`matzhn@nus.edu.sg`).

[‡]Department of Applied Mathematics, The Hong Kong Polytechnic University, Hung Hom, Hong Kong (`defeng.sun@polyu.edu.hk`). On leave from Department of Mathematics, National University of Singapore.

[§]Department of Mathematics, and Institute of Operations Research and Analytics, National University of Singapore, 10 Lower Kent Ridge Road, Singapore 119076 (`mattohkc@nus.edu.sg`).

# 1  Introduction

In this paper, we aim to design a fast algorithm for solving the following sparse group Lasso (SGLasso) problem:

$$\min_{x \in \mathbb{R}^n} \ \frac{1}{2}\|\mathcal{A}x - b\|^2 + \lambda_1 \|x\|_1 + \lambda_2 \sum_{l=1}^{g} w_l \|x_{G_l}\|, \tag{1}$$

where $\mathcal{A} : \mathbb{R}^n \to \mathbb{R}^m$ is a linear map, $b \in \mathbb{R}^m$ is the given response vector, $\lambda_1 \geq 0$ and $\lambda_2 \geq 0$ are regularization parameters. For $l = 1, 2, \ldots, g$, $w_l > 0$, and the set $G_l \subseteq \{1, 2, \ldots, n\}$ contains the indices corresponding to the $l$-th group of features. We denote the restriction of the vector $x$ to the index set $G_l$ as $x_{G_l}$. Here $\|\cdot\|$ and $\|\cdot\|_1$ denote the $\ell_2$ norm and $\ell_1$ norm, respectively. For convenience, we denote the SGLasso regularizer by the proper closed convex function $p(x) := \lambda_1 \|x\|_1 + \lambda_2 \sum_{l=1}^{g} w_l \|x_{G_l}\|$, $\forall\, x \in \mathbb{R}^n$.

In recent decades, high dimensional feature selection problems have become increasingly important, and the penalized regression models have been proven to be particularly useful for these feature selection problems. For many such problems in real applications, the number of predictors $n$ is much larger than the number of observations $m$. A notable example of the penalized regression model is the Lasso model that was first proposed by Tibshirani [41]. The problem (1) contains the Lasso problem as a special case if we take the parameter $\lambda_2 = 0$. On the other hand, by assuming that some prior information about the group structure of the underlying solution $x$ is known, Yuan and Lin [45] proposed the group Lasso model, i.e., in problem (1) with parameters $\lambda_1 = 0$ and $\lambda_2 \geq 0$. The group Lasso can select a small set of groups. However, it does not ensure sparsity within each group. For the purpose of achieving sparsity of groups and within each group, Friedman et al. [14] proposed the SGLasso model (1), potentially with overlaps between groups. Apart from the above penalized regression models, there exist a number of variants with different regularizers, such as the fused Lasso [42] and the network Lasso [17].

The SGLasso model has been widely applied to different fields, such as text processing, bioinformatics, signal interpretation, and object tracking (e.g., [20, 19, 32, 21, 10, 48]). Its wide ranging applications have inspired many researchers to design various algorithms for solving the SGLasso problem. These algorithms include the (accelerated) proximal gradient method (see e.g., [1, 45]), (randomized) block coordinate descent algorithm (see e.g., [35, 38, 34]), and alternating direction method of multipliers (see e.g., [2]). To the best of our knowledge, these existing algorithms are first order methods that are applied directly to the primal problem (1) and they hardly utilize any second order information. In contrast, we aim to design an efficient second order information based algorithm for solving the dual problem of the SGLasso problem (1). The motivations for adopting the dual approach will be explained in the following paragraphs.

Solving the SGLasso problem is especially challenging when there are overlapping groups because of the complex structure of the SGLasso regularizer $p$. The complicated composite structure of $p$ generally makes it impossible to compute its proximal mapping analytically. However, the efficient computation of such a proximal mapping is indispensable to a number of algorithms, and many of the papers mentioned in the last paragraph thus considered the simpler case of the non-overlapping SGLasso problem. As a first attempt to design a Hessian based algorithm for the SGLasso problem, we will also focus on the simpler case of the non-overlapping SGLasso problem. The non-overlapping case can be treated as a preliminary study towards the final goal of designing a Hessian based algorithm of solving the overlapping SGLasso problem. For the rest of this paper, we make the following blanket assumption.

**Assumption 1.1.** *The different groups $G_l$, $l = 1, 2, \ldots, g$ form a partition of $\{1, 2, \ldots, n\}$, i.e., $G_i \cap G_j = \emptyset$ for all $1 \leq i < j \leq g$, and $\cup_{l=1}^{g} G_l = \{1, 2, \ldots, n\}$.*

In order to solve the non-overlapping SGLasso problem, we aim to use the semismooth Newton (SSN) augmented Lagrangian (SSNAL) framework for solving the dual problem of (1). This approach is motivated by the superior numerical performance of the SSNAL when applied to the dual of the Lasso problem [24] and that of the fused Lasso problem [25]. We note that the objective functions of the Lasso and fused Lasso problems are piecewise linear-quadratic, and therefore as proven in [24, 25], both the primal and dual iterates generated by the augmented Lagrangian method (ALM) are asymptotically superlinearly convergent. It is this attractive convergence property that leads to the impressive numerical performance of the SSNAL. However, the regularizer $p$ in the objective function of the SGLasso problem (1) is no longer a polyhedral function due to the presence of the $\ell_2$ norm. As a result, the asymptotic superlinear convergence of both the primal and dual iterative sequences generated by the ALM are no longer guaranteed to hold by the existing theoretical results. Fortunately, by leveraging on the recent advances made in Cui, Sun, and Toh [8] on the analysis of the asymptotic R-superlinear convergence of the ALM for convex composite conic programming, we are able to establish the global linear convergence (with an arbitrary rate) of the primal iterative sequence, the dual infeasibility, and the dual function values generated by the ALM for the SGLasso problem. With this convergence result, we could expect the ALM to be highly efficient for solving the SGLasso problem.

The remaining challenge of designing an efficient ALM to solve (1) is in solving the subproblem in each iteration. As inspired by the success in [24, 25], we will design a highly efficient SSN method for solving the subproblem in each ALM iteration. The effectiveness of the SSN method relies greatly on the efficient computation of the generalized Jacobian of the proximal mapping associated with the SGLasso regularizer $p$. Thus a major contribution of this paper is to analyse the structure of the generalized Jacobian and its efficient computation. As far as we know, the elements in the generalized Jacobian of the proximal mapping of $p$ have not been derived before, and this paper aims to derive an explicit formula for them. We note that the SGLasso regularizer $p$ enjoys the "prox-decomposition" property [43], similar to the fused Lasso regularizer (see [25]). With the "prox-decomposition" property and some necessary properties for the $\ell_1$ norm and $\ell_2$ norm, we are able to derive an explicit formula for the generalized Jacobian of the proximal mapping of $p$. Based on the structure of the generalized Jacobian of the proximal mapping of $p$, we can derive a certain structured sparsity (which we name as the second order sparsity) of the Hessians associated with the objective function in each ALM subproblem to implement the SSN method efficiently. Moreover, the SSN method will be proven to have superlinear/quadratic convergence. In a nutshell, the globally fast linear convergence (with an arbitrary linear rate) of the ALM and the superlinear/quadratic convergence of the SSN method for solving each ALM subproblem can guarantee that our SSNAL is highly efficient and robust for solving large-scale SGLasso problems.

The rest of this paper is organized as follows. Section 2 demonstrates the decomposition property of the SGLasso regularizer and provides theoretical conditions for ensuring the global fast linear convergence of the ALM. The explicit formulation of the generalized Jacobian of the proximal mapping of the SGLasso regularizer is derived in section 3. In section 4, we design the semismooth Newton based augmented Lagrangian method (SSNAL) for solving the dual of the SGLasso problem (1) and derive our main convergence results. We will also present efficient techniques for implementing SSNAL. Section 5 evaluates the performance of SSNAL on both the synthetic and real data sets. Finally, concluding remarks are given in section 6.

**Notation.** For a linear map $\mathcal{A} : \mathbb{R}^n \to \mathbb{R}^m$, we denote its adjoint by $\mathcal{A}^*$. For any convex function $p$, we denote its conjugate function by $p^*$, i.e., $p^*(x) = \sup_z\{\langle x, z\rangle - p(z)\}$. For each $l \in \{1, 2, \ldots, g\}$, we define the linear operator $\mathcal{P}_l : \mathbb{R}^n \to \mathbb{R}^{|G_l|}$ by $\mathcal{P}_l x = x_{G_l}$. Let $s := \sum_{l=1}^g |G_l|$. Define $\mathcal{P} := [\mathcal{P}_1; \mathcal{P}_2; \ldots; \mathcal{P}_g] : \mathbb{R}^n \to \mathbb{R}^s$ and $\mathcal{B}_2 := \mathcal{B}_2^{\lambda_{2,1}} \times \cdots \times \mathcal{B}_2^{\lambda_{2,g}}$, where $\mathcal{B}_2^{\lambda_{2,l}} := \{u_l \in \mathbb{R}^{|G_l|} \mid \|u_l\| \leq \lambda_{2,l}\}$ and $\lambda_{2,l} := \lambda_2 w_l$. For a given closed convex set $\Omega$ and a vector $x$, denote the distance of $x$ to $\Omega$ by $\text{dist}(x, \Omega) := \inf_{x'\in\Omega}\{\|x - x'\|\}$ and the Euclidean projection of $x$ onto $\Omega$ by $\Pi_\Omega(x) := \arg\min_{x'\in\Omega}\{\|x - x'\|\}$. We define $\text{sign}(\cdot)$ in a component-wise fashion such that $\text{sign}(t) = 1$ if $t > 0$, $\text{sign}(t) = 0$ if $t = 0$, and $\text{sign}(t) = -1$ if $t < 0$. For any functions $f$ and $g$, define $(f \circ g)(\cdot) := f(g(\cdot))$. We denote the Hadamard product by $\odot$. For a given vector $x$, $\text{supp}(x)$ denotes the support of $x$, i.e., the set of indices such that $x_i \neq 0$. We denote the vector of all ones by $e$. For a matrix $A$ and a vector $a$, we denote by $\text{diag}(A)$ and $\text{Diag}(a)$ the diagonal vector of $A$ and the diagonal matrix whose diagonal elements are the components of $a$, respectively.

## 2 Preliminaries

In this section, we establish the decomposition property of the SGLasso regularizer $p$ and present some general error bound results. The SGLasso problem (1) can be written equivalently as follows:

$$(\text{P}) \quad \min_{x\in\mathbb{R}^n} \ h(x) := f(x) + p(x),$$

where $f(x) := \frac{1}{2}\|\mathcal{A}x-b\|^2$, $p(x) := \varphi(x)+\phi(x)$, $\varphi(x) := \lambda_1\|x\|_1$, and $\phi(x) := \sum_{l=1}^{g} \lambda_{2,l}\|x_{G_l}\|$ with $\lambda_{2,l} := \lambda_2 w_l$, $l = 1, 2, \ldots, g$. The dual problem of (P) takes the following form:

$$\text{(D)} \quad \begin{aligned} \max \quad & g(y,z) := -\langle b, y\rangle - \tfrac{1}{2}\|y\|^2 - p^*(z) \\ \text{s.t.} \quad & \mathcal{A}^*y + z = 0. \end{aligned}$$

In additon, the Karush-Kuhn-Tucker (KKT) optimality system associated with (P) and (D) is given by

$$\mathcal{A}x - y - b = 0,\ \operatorname{Prox}_p(x+z) - x = 0,\ \mathcal{A}^*y + z = 0, \qquad (2)$$

where the proximal mapping of $p$ is defined by:

$$\operatorname{Prox}_p(u) := \arg\min_x \left\{ p(x) + \frac{1}{2}\|x-u\|^2 \right\}, \ \forall\, u \in \mathbb{R}^n. \qquad (3)$$

For any given parameter $t > 0$ and a closed proper convex function $h$, the following Moreau identity will be frequently used:

$$\operatorname{Prox}_{th}(u) + t\operatorname{Prox}_{h^*/t}(u/t) = u. \qquad (4)$$

It is well known that the proximal mappings of $\ell_1$ norm and $\ell_2$ norm can be expressed as follows: for any given $c > 0$,

$$\operatorname{Prox}_{c\|\cdot\|_1}(u) = \operatorname{sign}(u) \odot \max\{|u| - ce, 0\},$$

$$\operatorname{Prox}_{c\|\cdot\|}(u) = \begin{cases} \frac{u}{\|u\|}\max\{\|u\|-c, 0\}, & \text{if } u \neq 0, \\ 0, & \text{otherwise}. \end{cases}$$

We will show in the next lemma that the proximal mappings of the $\ell_1$ norm and $\ell_2$ norm are strongly semismooth. The definition of the semismoothness was first introduced by Mifflin [31] for functionals and extended to vector-valued functions in [33].

**Lemma 2.1.** *For any $c > 0$, the proximal mappings $\operatorname{Prox}_{c\|\cdot\|_1}(\cdot)$ and $\operatorname{Prox}_{c\|\cdot\|}(\cdot)$ are strongly semismooth.*

*Proof.* Since $\operatorname{Prox}_{c\|\cdot\|_1}(\cdot)$ is a Lipschitz continuous piecewise affine function, it follows from [11, Proposition 7.4.7] that $\operatorname{Prox}_{c\|\cdot\|_1}(\cdot)$ is strongly semismooth everywhere. Next, we focus on the proximal mapping $\operatorname{Prox}_{c\|\cdot\|}(\cdot)$. From the definition of $\operatorname{Prox}_{c\|\cdot\|}(\cdot)$ and the fact that the projection of any vector onto the second order cone, i.e., the epigraph of the $\ell_2$ norm function, is strongly semismooth [5, Proposition 4.3], we can obtain the conclusion directly from [30, Theorem 4]. □

Next, we analyse the vital decomposition property, which is termed as "prox-decomposition" in [43], of the SGLasso regularizer $p$. In the next proposition, we show that the proximal mapping $\mathrm{Prox}_p(\cdot)$ of $p = \varphi + \phi$ can be decomposed into the composition of the proximal mappings $\mathrm{Prox}_\varphi(\cdot)$ and $\mathrm{Prox}_\phi(\cdot)$. With this decomposition property, we are able to compute $\mathrm{Prox}_p(\cdot)$ in a closed form. This decomposition result was proved in [44, Theorem 1], which is mainly an extension of that for the fused Lasso regularizer in [13]. Here, we give another short proof based on the systematic investigation in [43].

**Proposition 2.1.** *Under Assumption 1.1, it holds that*

$$\mathrm{Prox}_p(u) = \mathrm{Prox}_\phi \circ \mathrm{Prox}_\varphi(u), \ \forall\, u \in \mathbb{R}^n.$$

*Proof.* Under Assumption 1.1, the function $p$ has a separable structure. Hence, the problem (3) is separable for each group. Therefore, it is sufficient to prove that

$$\mathrm{Prox}_{\lambda_1\|\cdot\|_1 + \lambda_{2,l}\|\cdot\|}(u_l) = \mathrm{Prox}_{\lambda_{2,l}\|\cdot\|} \circ \mathrm{Prox}_{\lambda_1\|\cdot\|_1}(u_l), \ \forall\, u_l \in \mathbb{R}^{|G_l|}, \ l = 1, 2, \ldots, g.$$

By [43, Theorem 1], for each $l \in \{1, 2, \ldots, g\}$, it suffices to show that

$$\partial(\lambda_1 \|u_l\|_1) \subseteq \partial(\lambda_1 \|v_l\|_1), \ v_l := \mathrm{Prox}_{\lambda_{2,l}\|\cdot\|}(u_l), \ \forall\, u_l \in \mathbb{R}^{|G_l|}.$$

For any given $u_l \in \mathbb{R}^{|G_l|}$, we discuss the following two cases.

Case 1: If $\|u_l\| \leq \lambda_{2,l}$, then $v_l = 0$. It follows that $\partial(\lambda_1 \|v_l\|_1) = [-\lambda_1, \lambda_1]^{|G_l|}$, which obviously contains $\partial(\lambda_1 \|u_l\|_1)$.

Case 2: If $\|u_l\| > \lambda_{2,l}$, then $v_l = (1 - \lambda_{2,l}/\|u_l\|)u_l$, which implies that $\mathrm{sign}(v_l) = \mathrm{sign}(u_l)$. Thus, it holds that $\partial(\lambda_1 \|u_l\|_1) = \partial(\lambda_1 \|v_l\|_1)$.

Hence, the proof is completed. $\square$

Consider an arbitrary point $u \in \mathbb{R}^n$. Based on the above proposition, we are now ready to compute $\mathrm{Prox}_p(u)$ explicitly. Let $v := \mathrm{Prox}_\varphi(u)$. For each group $G_l$, $l = 1, 2, \ldots, g$, it holds that

$$\arg\min_{x_{G_l}} \left\{ \lambda_{2,l} \|x_{G_l}\| + \frac{1}{2}\|x_{G_l} - v_{G_l}\|^2 \right\} = v_{G_l} - \Pi_{\mathcal{B}_2^{\lambda_{2,l}}}(v_{G_l}).$$

That is, $\mathcal{P}_l \mathrm{Prox}_\phi(v) = \mathcal{P}_l v - \Pi_{\mathcal{B}_2^{\lambda_{2,l}}}(\mathcal{P}_l v)$. Therefore, we have

$$\mathrm{Prox}_p(u) = \mathrm{Prox}_\phi(v) = v - \mathcal{P}^* \Pi_{\mathcal{B}_2}(\mathcal{P}v). \qquad (5)$$

For the rest of this section, we introduce some error bound results that will be used later in the convergence rate analysis. Define the proximal residual function $\mathcal{R}: \mathbb{R}^n \to \mathbb{R}^n$ by

$$\mathcal{R}(x) := x - \mathrm{Prox}_p(x - \nabla f(x)), \quad \forall\, x \in \mathbb{R}^n \qquad (6)$$

and the set-valued map $\mathcal{T}: \mathbb{R}^n \rightrightarrows \mathbb{R}^n$ by

$$\mathcal{T}(x) := \{v \in \mathbb{R}^n \,|\, v \in \nabla f(x) + \partial p(x)\}, \quad \forall\, x \in \mathbb{R}^n. \tag{7}$$

Suppose that $\lambda_1 + \lambda_2 > 0$. Let $\Omega_P$ be the optimal solution set of $(P)$. Since $f$ is nonnegative on $\mathbb{R}^n$, it is easy to obtain that $h(x) \to +\infty$ as $\|x\| \to +\infty$. Thus, $\Omega_P$ is a compact convex set. The first order optimality condition of $(P)$ implies that $\bar{x} \in \Omega_P$ is equivalent to $0 \in \mathcal{T}(\bar{x})$, which in turn is equivalent to $\mathcal{R}(\bar{x}) = 0$. It is proved in [46, Theorem 1] that the local error bound condition (in the sense of Luo and Tseng [28]) holds around the optimal set $\Omega_P$, i.e., for every $\xi \geq \inf_x h(x)$, there exist positive scalars $\kappa_0$ and $\delta_0$ such that

$$\mathrm{dist}(x, \Omega_P) \leq \kappa_0 \|\mathcal{R}(x)\|, \ \forall\, x \in \mathbb{R}^n \text{ satisfying } h(x) \leq \xi \text{ and } \|\mathcal{R}(x)\| \leq \delta_0. \tag{8}$$

Therefore, by using the facts that $\Omega_P$ is compact and that $\mathcal{R}$ is continuous, we know that for any $r_1 > 0$, there exists $\kappa_1 > 0$ such that

$$\mathrm{dist}(x, \Omega_P) \leq \kappa_1 \|\mathcal{R}(x)\|, \ \forall\, x \in \mathbb{R}^n \text{ satisfying } \mathrm{dist}\,(x, \Omega_P) \leq r_1. \tag{9}$$

Furthermore, by mimicking the proofs in [9, Theorem 3.1] or [7, Proposition 2.4] and noting that $\Omega_P$ is a compact set, we can obtain the following result with no difficulty.

**Proposition 2.2.** *For any $r > 0$, there exists $\kappa > 0$ such that*

$$\mathrm{dist}(x, \Omega_P) \leq \kappa \,\mathrm{dist}(0, \mathcal{T}(x)), \ \forall\, x \in \mathbb{R}^n \text{ satisfying } \mathrm{dist}\,(x, \Omega_P) \leq r.$$

## 3 Generalized Jacobian of $\mathrm{Prox}_p(\cdot)$

In this section, we shall analyse the generalized Jacobian of the proximal mapping $\mathrm{Prox}_p(\cdot)$ of the SGLasso regularizer $p$. From Proposition 2.1, for any $u \in \mathbb{R}^n$, we have

$$\mathrm{Prox}_p(u) = \mathrm{Prox}_\phi(\mathrm{Prox}_\varphi(u)).$$

At the first glance, we may try to apply the chain rule in deriving the generalized Jacobian of $\mathrm{Prox}_p(\cdot)$. Indeed it was illustrated in [39] that under certain conditions, the generalized Jacobian for composite functions can be obtained by the chain rule in a similar fashion as in finding the ordinary Jacobian for composite smooth functions. Specifically, if the conditions in [39, Lemma 2.1] hold, then we could have obtained by the chain rule the following B-subdifferential [6], which is a subset of the generalized Jacobian,

$$\partial_B \mathrm{Prox}_p(u) = \left\{ \tilde{\Theta} \cdot \Theta \,\middle|\, \tilde{\Theta} \in \partial_B \mathrm{Prox}_\phi(v),\ \Theta \in \partial_B \mathrm{Prox}_\varphi(u),\ v = \mathrm{Prox}_\varphi(u) \right\}.$$

However, the conditions in [39, Lemma 2.1] may not hold in our context, and consequently the above equation is usually invalid. Therefore, the generalized Jacobian of $\mathrm{Prox}_p(\cdot)$ is nontrivial to obtain, and we have to find an alternative surrogate to bypass this difficulty. The challenge just highlighted also appeared in [25] when analysing the generalized Jacobian of the proximal mapping of the fused Lasso regularizer. In that work, a general definition (see [25, Definition 1]) of "semismoothness with respect to a multifunction" was adopted, and such a multifunction was constructed to play the role of the generalized Jacobian. Here, we shall use the same strategy, and our task now is to identify such a multifunction.

Before characterizing the multifunction relating to the semismoothness, based on the fact in (5) that $\text{Prox}_\phi(v) = v - \mathcal{P}^*\Pi_{\mathcal{B}_2}(\mathcal{P}v)$, $\forall v \in \mathbb{R}^n$, we define the following alternative for the generalized Jacobian of $\text{Prox}_\phi(\cdot)$:

$$\hat{\partial}\text{Prox}_\phi(v) := \left\{ I - \mathcal{P}^*\Sigma\mathcal{P} \,\big|\, \Sigma = \text{Diag}(\Sigma_1, \ldots, \Sigma_g), \Sigma_l \in \partial\Pi_{\mathcal{B}_2^{\lambda_{2,l}}}(v_{G_l}), l = 1, 2, \ldots, g \right\}.$$

It can be observed that the main part of $\hat{\partial}\text{Prox}_\phi(\cdot)$ is the block diagonal matrix $\Sigma$, of which each block is the generalized Jacobian of a projection operator onto an $\ell_2$-norm ball. Since $\partial\Pi_{\mathcal{B}_2^{\lambda_{2,l}}}(\cdot)$ admits a closed form expression, so does $\hat{\partial}\text{Prox}_\phi(\cdot)$. Now, we are in a position to present the following multifunction $\mathcal{M} : \mathbb{R}^n \rightrightarrows \mathbb{R}^{n \times n}$ and regard it as the surrogate generalized Jacobian of $\text{Prox}_p(\cdot)$ at any $u \in \mathbb{R}^n$:

$$\mathcal{M}(u) := \left\{ (I - \mathcal{P}^*\Sigma\mathcal{P})\Theta \,\bigg|\, \begin{array}{l} \Sigma = \text{Diag}(\Sigma_1, \ldots, \Sigma_g), \Sigma_l \in \partial\Pi_{\mathcal{B}_2^{\lambda_{2,l}}}(v_{G_l}), l = 1, 2, \ldots, g, \\ v = \text{Prox}_\varphi(u), \Theta \in \partial\text{Prox}_\varphi(u) \end{array} \right\}. \tag{10}$$

**Remark 3.1.** *In numerical computations, for any $u \in \mathbb{R}^n$, one needs to construct at least one element in $\mathcal{M}(u)$ explicitly. For $l = 1, 2, \ldots, g$ and $v_l \in \mathbb{R}^{|G_l|}$, the projection onto an $\ell_2$-norm ball and its generalized Jacobian are given as follows, respectively:*

$$\Pi_{\mathcal{B}_2^{\lambda_{2,l}}}(v_l) = \begin{cases} \lambda_{2,l}\frac{v_l}{\|v_l\|}, & \text{if } \|v_l\| > \lambda_{2,l}, \\ v_l, & \text{otherwise}, \end{cases} \tag{11}$$

$$\partial\Pi_{\mathcal{B}_2^{\lambda_{2,l}}}(v_l) = \begin{cases} \frac{\lambda_{2,l}}{\|v_l\|}(I - \frac{v_l v_l^T}{\|v_l\|^2}), & \text{if } \|v_l\| > \lambda_{2,l}, \\ \{I - t\frac{v_l v_l^T}{(\lambda_{2,l})^2} \,|\, 0 \leq t \leq 1\}, & \text{if } \|v_l\| = \lambda_{2,l}, \\ I, & \text{if } \|v_l\| < \lambda_{2,l}. \end{cases} \tag{12}$$

*Define a vector $\theta \in \mathbb{R}^n$ and construct a matrix $\Theta \in \partial\text{Prox}_\varphi(u)$ as follow:*

$$\theta_i = \begin{cases} 0, & \text{if } |u_i| \leq \lambda_1, \\ 1, & \text{otherwise}, \end{cases} \quad i = 1, \ldots, n. \tag{13}$$

*Then it is obvious that*

$$\Theta = \text{Diag}(\theta) \in \partial\text{Prox}_\varphi(u). \tag{14}$$

The following main theorem of this section justifies why $\mathcal{M}(u)$ in (10) can be treated as the surrogate generalized Jacobian of $\text{Prox}_p(\cdot)$ at $u$. That is, it shows that the proximal mapping $\text{Prox}_p$ is strongly semismooth on $\mathbb{R}^n$ with respect to the multifunction $\mathcal{M}$ defined in (10) (cf. Definition 1 in [25]).

**Theorem 3.1.** *Assume that Assumption 1.1 holds. Let $u \in \mathbb{R}^n$. Then the multifunction $\mathcal{M}$, defined in (10), is a nonempty compact valued upper-semicontinuous multifunction, and for any $M \in \mathcal{M}(u)$, $M$ is symmetric and positive semidefinite. Moreover, for any $M \in \mathcal{M}(w)$ with $w \to u$,*

$$\text{Prox}_p(w) - \text{Prox}_p(u) - M(w - u) = O(\|w - u\|^2). \tag{15}$$

*Proof.* By Lemma 2.1, Proposition 2.1 and [11, Theorem 7.5.17], one can deduce that the point-to-set map $\mathcal{M}$ has nonempty compact images and is upper-semicontinuous, and equation (15) holds. It remains to show that $M$ is symmetric and positive semidefinite for any $M \in \mathcal{M}(u)$. Denote $v := \text{Prox}_\varphi(u)$. Take $M \in \mathcal{M}(u)$ arbitrarily. Then, there exist $\Sigma_l \in \partial \Pi_{\mathcal{B}_2^{\lambda_{2,l}}}(v_{G_l})$, $l = 1, 2, \ldots, g$ and $\Theta = \text{Diag}(\theta) \in \partial \text{Prox}_\varphi(u)$, defined in (12) and (14), respectively, such that

$$M = \sum_{l=1}^{g} \mathcal{P}_l^*(I - \Sigma_l)\mathcal{P}_l\Theta.$$

It suffices to show that $\mathcal{P}_l^*(I - \Sigma_l)\mathcal{P}_l\Theta$ is symmetric and positive semidefinite for any $l \in \{1, 2, \ldots, g\}$. Denote the index sets

$$\Xi_l := G_l \cap \text{supp}(v) = \{i \in G_l \mid \theta_i = 1\},\ l = 1, 2, \ldots, g. \tag{16}$$

For simplicity, we write $v_{G_l}$ as $v_l$ in the following proof.

Case 1: $\|v_l\| < \lambda_{2,l}$. By (12), $I - \Sigma_l = 0$.

Case 2: $\|v_l\| = \lambda_{2,l}$. By (12), there exists some $t \in [0, 1]$ such that

$$\mathcal{P}_l^*(I - \Sigma_l)\mathcal{P}_l\Theta = \frac{t}{(\lambda_{2,l})^2}(\mathcal{P}_l^* v_l)(\mathcal{P}_l^* v_l)^T \Theta.$$

By the definition of $\mathcal{P}_l$, we deduce that $\text{supp}(\mathcal{P}_l^* v_l) \subseteq \Xi_l$. It follows from (16) that $(\mathcal{P}_l^* v_l)^T \Theta = (\mathcal{P}_l^* v_l)^T$. That is,

$$\mathcal{P}_l^*(I - \Sigma_l)\mathcal{P}_l\Theta = \frac{t}{(\lambda_{2,l})^2}(\mathcal{P}_l^* v_l)(\mathcal{P}_l^* v_l)^T,$$

which is symmetric and positive semidefinite.

Case 3: $\|v_l\| > \lambda_{2,l}$. From (12) and the proof in case 2, we have

$$\mathcal{P}_l^*(I - \Sigma_l)\mathcal{P}_l\Theta = \mathcal{P}_l^*\Big(I - \frac{\lambda_{2,l}}{\|v_l\|}\big(I - \frac{v_l v_l^T}{\|v_l\|^2}\big)\Big)\mathcal{P}_l\Theta$$

$$= \Big(1 - \frac{\lambda_{2,l}}{\|v_l\|}\Big)\mathcal{P}_l^*\mathcal{P}_l\Theta + \frac{\lambda_{2,l}}{\|v_l\|^3}(\mathcal{P}_l^* v_l)(\mathcal{P}_l^* v_l)^T\Theta$$

$$= \Big(1 - \frac{\lambda_{2,l}}{\|v_l\|}\Big)\mathcal{P}_l^*\mathcal{P}_l\Theta + \frac{\lambda_{2,l}}{\|v_l\|^3}(\mathcal{P}_l^* v_l)(\mathcal{P}_l^* v_l)^T.$$

Since both $\mathcal{P}_l^*\mathcal{P}_l$ and $\Theta$ are diagonal, it holds that $\mathcal{P}_l^*(I - \Sigma_l)\mathcal{P}_l\Theta$ is symmetric. Furthermore, it is obvious that $\mathcal{P}_l^*\mathcal{P}_l\Theta$ is positive semidefinite. Therefore, the last equality implies that $\mathcal{P}_l^*(I - \Sigma_l)\mathcal{P}_l\Theta$ is positive semidefinite.

In summary, we have shown that $M$ is symmetric and positive semidefinite. $\square$

## 4 An inexact semismooth Newton based augmented Lagrangian method

In this section, we shall design an inexact semismooth Newton based augmented Lagrangian method for solving problem (D), the dual of the SGLasso problem (1). Here we always assume that $\lambda_1 + \lambda_2 > 0$. Write (D) equivalently in the following

$$\begin{aligned} \min \quad & \langle b, y \rangle + \tfrac{1}{2}\|y\|^2 + p^*(z) \\ \text{s.t.} \quad & \mathcal{A}^* y + z = 0. \end{aligned} \qquad (17)$$

For $\sigma > 0$, the augmented Lagrangian function associated with (17) is given by

$$\mathcal{L}_\sigma(y, z; x) = \langle b, y \rangle + \frac{1}{2}\|y\|^2 + p^*(z) + \frac{\sigma}{2}\|\mathcal{A}^* y + z - \sigma^{-1}x\|^2 - \frac{1}{2\sigma}\|x\|^2. \qquad (18)$$

The $k$-th iteration of the augmented Lagrangian method is given as follows:

$$\begin{cases} (y^{k+1}, z^{k+1}) \approx \arg\min_{y,z}\{\mathcal{L}_{\sigma_k}(y, z; x^k)\}, \\ x^{k+1} = x^k - \sigma_k(\mathcal{A}^* y^{k+1} + z^{k+1}),\ k \geq 0. \end{cases}$$

In each iteration, the most expensive step is to solve the following subproblem:

$$\min_{y,z}\{\mathcal{L}_{\sigma_k}(y, z; x^k)\}. \qquad (19)$$

Obviously, the subproblem (19) admits a unique optimal solution. For any $y \in \mathbb{R}^m$, define

$$\psi_k(y) := \inf_z \mathcal{L}_{\sigma_k}(y, z; x^k)$$

$$= \langle b, y \rangle + \frac{1}{2}\|y\|^2 + p^*(\mathrm{Prox}_{p^*/\sigma_k}(\sigma_k^{-1}x^k - \mathcal{A}^* y)) + \frac{\sigma_k}{2}\|\mathrm{Prox}_p(\sigma_k^{-1}x^k - \mathcal{A}^* y)\|^2$$

$$- \frac{1}{2\sigma_k}\|x^k\|^2. \qquad (20)$$

Then, $(y^{k+1}, z^{k+1}) \approx \arg\min_{y,z}\{\mathcal{L}_{\sigma_k}(y, z; x^k)\}$ can be computed as follows:

$$y^{k+1} \approx \arg\min_y \psi_k(y) \quad \text{and} \quad z^{k+1} = \text{Prox}_{p^*/\sigma_k}(\sigma_k^{-1}x^k - \mathcal{A}^* y^{k+1}). \tag{21}$$

Now, we propose an inexact augmented Lagrangian method for solving (17).

---
**Algorithm 1:** An inexact augmented Lagrangian method for solving (17)

---
Let $\sigma_0 > 0$ be a given parameter. Choose $(y^0, z^0, x^0) \in \mathbb{R}^m \times \mathbb{R}^n \times \mathbb{R}^n$. Iterate the following steps for $k = 0, 1, \ldots$.

**Step 1.** Compute
$$(y^{k+1}, z^{k+1}) \approx \arg\min_{y,z}\{\mathcal{L}_{\sigma_k}(y, z; x^k)\} \tag{22}$$
via (21).

**Step 2.** Compute
$$x^{k+1} = x^k - \sigma_k(\mathcal{A}^* y^{k+1} + z^{k+1}) = \sigma_k \text{Prox}_p(\sigma_k^{-1} x^k - \mathcal{A}^* y^{k+1}). \tag{23}$$

**Step 3.** Update $\sigma_{k+1} \uparrow \sigma_\infty \leq \infty$.

---

Given nonnegative summable sequences $\{\varepsilon_k\}$ and $\{\delta_k\}$ such that $\delta_k < 1$ for all $k \geq 0$, we estimate the accuracy of the approximate solution $(y^{k+1}, z^{k+1})$ of (22) via the standard stopping criteria studied in [36]:

(A) $\mathcal{L}_{\sigma_k}(y^{k+1}, z^{k+1}; x^k) - \inf_{y,z} \mathcal{L}_{\sigma_k}(y, z; x^k) \leq \varepsilon_k^2/2\sigma_k$,

(B) $\mathcal{L}_{\sigma_k}(y^{k+1}, z^{k+1}; x^k) - \inf_{y,z} \mathcal{L}_{\sigma_k}(y, z; x^k) \leq (\delta_k^2/2\sigma_k)\|x^{k+1} - x^k\|^2$.

Since $\psi_k$ is strongly convex with modulus 1, one has the estimate

$$\mathcal{L}_{\sigma_k}(y^{k+1}, z^{k+1}; x^k) - \inf_{y,z} \mathcal{L}_{\sigma_k}(y, z; x^k) = \psi_k(y^{k+1}) - \inf \psi_k \leq \frac{1}{2}\|\nabla \psi_k(y^{k+1})\|^2.$$

Therefore, the above stopping criteria (A) and (B) can be replaced by the following easy-to-check criteria, respectively,

(A′) $\|\nabla \psi_k(y^{k+1})\| \leq \varepsilon_k/\sqrt{\sigma_k}$,

(B′) $\|\nabla \psi_k(y^{k+1})\| \leq (\delta_k/\sqrt{\sigma_k})\|x^{k+1} - x^k\|$.

## 4.1 Convergence rates for Algorithm 1

In this section, we shall analyse the linear convergence at an arbitrarily fast rate of the inexact augmented Lagrangian method for solving problem (17).

For the nonnegative summable sequence $\{\varepsilon_k\}$ in the stopping criterion $(A')$, we introduce a scalar $\alpha$ such that $\sum_{k=0}^{\infty} \varepsilon_k \leq \alpha$. Let $r$ be any given positive scalar satisfying $r > \alpha$. It follows from Proposition 2.2 that there exists a positive scalar $\kappa$ such that

$$\mathrm{dist}(x, \Omega_P) \leq \kappa \, \mathrm{dist}(0, \mathcal{T}(x)), \ \forall x \in \mathbb{R}^n \text{ satisfying } \mathrm{dist}(x, \Omega_P) \leq r. \tag{24}$$

The next lemma measures the distance of each primal iterate generated by Algorithm 1 to the optimal solution set $\Omega_P$. The proof of Lemma 4.1 is mainly based on [8, Proposition 1(c)], which itself is an extension of [37, Theorem 2] and [29, Theorem 2.1]. Compared to the proof in [8, Proposition 1(c)], the following lemma uses (24) instead of the calmness condition of $\mathcal{T}^{-1}$ at the origin for some $\bar{x} \in \Omega_P$.

**Lemma 4.1.** *Suppose that the initial point $x^0 \in \mathbb{R}^n$ satisfies $\mathrm{dist}(x^0, \Omega_P) \leq r - \alpha$. Let $\{x^k\}$ be any infinite sequence generated by Algorithm 1 under criteria $(A')$ and $(B')$. Then for all $k \geq 0$, one has*

$$\mathrm{dist}(x^{k+1}, \Omega_P) \leq \mu_k \, \mathrm{dist}(x^k, \Omega_P),$$

*where $\mu_k := [\delta_k + (1 + \delta_k)\kappa/\sqrt{\kappa^2 + \sigma_k^2}]/(1 - \delta_k)$ and $\kappa$ is from (24).*

*Proof.* Denote the proximal point mapping by $P_k := (\mathcal{I} + \sigma_k \mathcal{T})^{-1}$. Criterion $(A')$, [36, Proposition 6], and [37, Equation (2.6)] imply that

$$\|x^k - \Pi_{\Omega_P}(x^0)\| \leq \|x^0 - \Pi_{\Omega_P}(x^0)\| + \sum_{i=0}^{k-1} \varepsilon_i \leq \|x^0 - \Pi_{\Omega_P}(x^0)\| + \alpha, \ \forall k \geq 0.$$

Consequently, $\mathrm{dist}(x^k, \Omega_P) \leq \mathrm{dist}(x^0, \Omega_P) + \alpha \leq r, \ \forall k \geq 0$. Moreover, one has

$$\|P_k(x^k) - \Pi_{\Omega_P}(x^k)\| = \|P_k(x^k) - P_k(\Pi_{\Omega_P}(x^k))\| \leq \|x^k - \Pi_{\Omega_P}(x^k)\| \leq r,$$

which implies that

$$\mathrm{dist}(P_k(x^k), \Omega_P) \leq r, \ \forall k \geq 0.$$

Additionally, it is shown in [37, Proposition 1(a)] that

$$P_k(x^k) \in \mathcal{T}^{-1}((x^k - P_k(x^k))/\sigma_k), \ \forall k \geq 0.$$

Then it follows from (24) that

$$\mathrm{dist}(P_k(x^k), \Omega_P) \leq \kappa \, \mathrm{dist}(0, \mathcal{T}(P_k(x^k))) \leq (\kappa/\sigma_k)\|x^k - P_k(x^k)\|, \ \forall k \geq 0.$$

Therefore, from the proof in [8, Proposition 1 (c)], for all $k \geq 0$, we obtain that

$$\mathrm{dist}(P_k(x^k), \Omega_P) \leq (\kappa/\sqrt{\kappa^2 + \sigma_k^2}) \, \mathrm{dist}(x^k, \Omega_P)$$

and that

$$\|x^{k+1} - \Pi_{\Omega_P}(P_k(x^k))\|$$
$$\leq \delta_k \|x^{k+1} - \Pi_{\Omega_P}(P_k(x^k))\| + \left(\delta_k + (1+\delta_k)\kappa/\sqrt{\kappa^2 + \sigma_k^2}\right) \operatorname{dist}(x^k, \Omega_P).$$

This, together with the fact that $\operatorname{dist}(x^{k+1}, \Omega_P) \leq \|x^{k+1} - \Pi_{\Omega_P}(P_k(x^k))\|$, $\forall\, k \geq 0$, completes the proof. □

While the global convergence of Algorihm 1 follows from [36, 29] directly, the conditions required in [36, 29] to guarantee the local linear convergence of both $\{x^k\}$ and $\{(y^k, z^k)\}$ may no longer hold for the SGLasso problem due to the non-polyhedral property of the $\ell_2$ norm function. Fortunately, the new results established in [8] on the convergence rates of the ALM allow us to establish the following theorem, which proves the global Q-linear convergence of the primal sequence $\{x^k\}$ and the global R-linear convergence of the dual infeasibility and the dual objective values. Furthermore, the linear rates can be arbitrarily fast if the penalty parameter $\sigma_k$ is chosen sufficiently large.

**Theorem 4.1.** *Let $\{(y^k, z^k, x^k)\}$ be an infinite sequence generated by Algorithm 1 under stopping criterion $(A')$. Then, the sequence $\{x^k\}$ converges to some $\bar{x} \in \Omega_P$, and the sequence $\{(y^k, z^k)\}$ converges to the unique optimal solution of $(D)$.*

*Furthermore, if criterion $(B')$ is also executed in Algorithm 1 and the initial point $x^0 \in \mathbb{R}^n$ satisfies $\operatorname{dist}(x^0, \Omega_P) \leq r - \alpha$, then for all $k \geq 0$, we have*

$$\operatorname{dist}(x^{k+1}, \Omega_P) \leq \mu_k \operatorname{dist}(x^k, \Omega_P), \tag{25a}$$
$$\|\mathcal{A}^* y^{k+1} + z^{k+1}\| \leq \mu_k' \operatorname{dist}(x^k, \Omega_P), \tag{25b}$$
$$\sup(D) - g(y^{k+1}, z^{k+1}) \leq \mu_k'' \operatorname{dist}(x^k, \Omega_P), \tag{25c}$$

*where*

$$\mu_k := \left[\delta_k + (1+\delta_k)\kappa/\sqrt{\kappa^2 + \sigma_k^2}\right]/(1-\delta_k),$$
$$\mu_k' := 1/[(1-\delta_k)\sigma_k],$$
$$\mu_k'' := [\delta_k^2 \|x^{k+1} - x^k\| + \|x^{k+1}\| + \|x^k\|]/[2(1-\delta_k)\sigma_k],$$

*and $\kappa$ is from (24). Moreover, $\mu_k$, $\mu_k'$, and $\mu_k''$ go to 0 if $\sigma_k \uparrow \sigma_\infty = +\infty$.*

*Proof.* The statements on the global convergence just follow from [36, Theorem 5] or [8, Proposition 2]. Inequality (25a) is a direct consequence of Lemma 4.1. Inequality (25c) can be obtained by [8, Proposition 2 (5b)]. From the updating formula (23) of $x^{k+1}$, we deduce that

$$\|\mathcal{A}^* y^{k+1} + z^{k+1}\| = \sigma_k^{-1} \|x^{k+1} - x^k\|,$$

which, together with [8, Lemma 3], implies that (25b) holds. This completes the proof. □

**Remark 4.1.** *Assume that all the conditions in Theorem 4.1 are satisfied. Since the primal objective function h is Lipschitz continuous on any compact set, there exists a constant $L > 0$ such that h is Lipschitz continuous on the set $\{x \in \mathbb{R}^n \,|\, \mathrm{dist}(x, \Omega_P) \leq r\}$ with modulus L. Therefore, one can obtain from Theorem 4.1 that for all $k \geq 0$,*

$$h(x^{k+1}) - \inf(P) \leq L\mathrm{dist}(x^{k+1}, \Omega_P) \leq L\mu_k \mathrm{dist}(x^k, \Omega_P).$$

*This inequality, together with (25c) and the strong duality theorem, implies that*

$$h(x^{k+1}) - g(y^{k+1}, z^{k+1}) \leq (L\mu_k + \mu_k'')\mathrm{dist}\,(x^k, \Omega_P),$$

*which means that the duality gap converges to zero R-linearly at an arbitrary linear rate if $\sigma_k$ is sufficiently large and R-superlinearly if $\sigma_k \uparrow \sigma_\infty = +\infty$.*

## 4.2 A semismooth Newton method for solving the subproblem (21)

In this subsection, we propose an efficient semismooth Newton (SSN) method for solving the subproblem (21). As already mentioned earlier, having an efficient method for solving (21) is critical to the efficiency of Algorithm 1. In each iteration, we have to solve the following problem, for any given $\sigma > 0$ and fixed $\tilde{x}$,

$$\min_y \left\{ \psi(y) := \langle b, y \rangle + \frac{1}{2}\|y\|^2 + p^*(\mathrm{Prox}_{p^*/\sigma}(\sigma^{-1}\tilde{x} - \mathcal{A}^*y)) + \frac{\sigma}{2}\|\mathrm{Prox}_p(\sigma^{-1}\tilde{x} - \mathcal{A}^*y)\|^2 \right\}. \tag{26}$$

Note that $\psi(\cdot)$ is strongly convex and continuously differentiable with

$$\nabla \psi(y) = b + y - \sigma \mathcal{A}\mathrm{Prox}_p(\sigma^{-1}\tilde{x} - \mathcal{A}^*y).$$

Thus, the unique solution $\bar{y}$ of (26) can be obtained by solving the following nonsmooth equation

$$\nabla \psi(y) = 0. \tag{27}$$

Generally, to solve

$$F(x) = 0,$$

where $F : \mathbb{R}^m \to \mathbb{R}^m$ is a locally Lipschitz continuous function, one can employ the following SSN method:

$$x^{k+1} = x^k - V_k^{-1} F(x^k),$$

where $V_k \in \partial F(x^k)$, and $\partial F(x^k)$ denotes the Clarke generalized Jacobian [6, Definition 2.6.1] of $F$ at $x^k$. For more details about the SSN method, we refer the reader to [22, 23, 33, 40, 47] and the references therein. In particular, existing studies such as [33, 40] used the Clarke generalized Jacobian $V_k \in \partial F(x^k)$ in the updating scheme and established correspondingly the convergence results of the SSN method.

We should point out again that characterizing $\partial(\nabla\psi)(\cdot)$ is a difficult task to accomplish. In section 3, we have constructed a multifunction $\mathcal{M}$, which is used as a surrogate of the generalized Jacobian $\partial\text{Prox}_p$. Besides, it is illustrated in Theorem 3.1 that $\text{Prox}_p$ is strongly semismooth with respect to the multifunction $\mathcal{M}$. Likewise, we define a multifunction $\mathcal{V}: \mathbb{R}^m \rightrightarrows \mathbb{R}^{m\times m}$ as follows:

$$\mathcal{V}(y) := \left\{V \,|\, V = I + \sigma\mathcal{A}M\mathcal{A}^*, \ M \in \mathcal{M}(\sigma^{-1}\tilde{x} - \mathcal{A}^*y)\right\},$$

where $\mathcal{M}(\cdot)$ is defined in (10). It follows from Theorem 3.1 and [11, Theorem 7.5.17] that (i) $\mathcal{V}$ is a nonempty compact valued upper-semicontinuous multifunction; (ii) $\nabla\psi$ is strongly semismooth on $\mathbb{R}^m$ with respect to the multifunction $\mathcal{V}$; (iii) every matrix in the set $\mathcal{V}(\cdot)$ is symmetric and positive definite.

With the above analysis, we are ready to design the following SSN method for solving (27).

---

**Algorithm 2:** A semismooth Newton method for solving (27)

Given $\mu \in (0, 1/2)$, $\bar{\eta} \in (0,1)$, $\tau \in (0,1]$, and $\delta \in (0,1)$. Choose $y^0 \in \mathbb{R}^m$. Iterate the following steps for $j = 0, 1, \ldots$.

**Step 1.** Choose $M_j \in \mathcal{M}(\sigma^{-1}\tilde{x} - \mathcal{A}^*y^j)$. Let $V_j = I + \sigma\mathcal{A}M_j\mathcal{A}^*$. Solve the following linear system

$$V_j d = -\nabla\psi(y^j) \tag{28}$$

exactly or by the conjugate gradient (CG) algorithm to find $d^j$ such that $\|V_j d^j + \nabla\psi(y^j)\| \leq \min(\bar{\eta}, \|\nabla\psi(y^j)\|^{1+\tau})$.

**Step 2.** (Line search) Set $\alpha_j = \delta^{m_j}$, where $m_j$ is the smallest nonnegative integer $m$ for which

$$\psi(y^j + \delta^m d^j) \leq \psi(y^j) + \mu\delta^m \langle \nabla\psi(y^j), d^j \rangle.$$

**Step 3.** Set $y^{j+1} = y^j + \alpha_j d^j$.

---

The following convergence theorem for Algorithm 2 can be obtained directly from [25, Theorem 3].

**Theorem 4.2.** *Let $\{y^j\}$ be the sequence generated by Algorithm 2. Then $\{y^j\}$ is well-defined and converges to the unique solution $\bar{y}$ of (26). Moreover, the convergence rate is at least superlinear:*

$$\|y^{j+1} - \bar{y}\| = O(\|y^j - \bar{y}\|^{1+\tau}),$$

*where $\tau \in (0,1]$ is the parameter given in Algorithm 2.*

## 4.3 Efficient techniques for solving the linear system (28)

In this section, we analyse the sparsity structure of the matrix in the linear system (28) and design sophisticated numerical techniques for solving the large-scale linear systems involved in the SSN method. These techniques were first applied in [24] which took full advantage of the second order sparsity of the underlying problem.

As can be seen, the most expensive step in each iteration of Algorithm 2 is in solving the linear system (28). Let $(\tilde{x}, y) \in \mathbb{R}^n \times \mathbb{R}^m$ and $\sigma > 0$ be given. The linear system (28) has the following form:
$$(I + \sigma A M A^T)d = -\nabla \psi(y), \qquad (29)$$
where $A$ denotes the matrix representation of the linear operator $\mathcal{A}$, and $M \in \mathcal{M}(u)$ with $u = \sigma^{-1}\tilde{x} - A^T y$. With the fact that $A$ is an $m$ by $n$ matrix and $M$ is an $n$ by $n$ matrix, the cost of naively computing $AMA^T$ is $O(mn(m+n))$. Similarly, for any vector $d \in \mathbb{R}^m$, the cost of naively computing the matrix-vector product $AMA^T d$ is $O(mn)$. Since the cost of naively computing the coefficient matrix $I + \sigma AMA^T$ and that of multiplying a vector by the coefficient matrix $I + \sigma AMA^T$ are excessively demanding, common linear system solvers, such as the Cholesky decomposition and the conjugate gradient method, will be extremely slow (if possible at all) in solving the linear system (29) arising from large-scale problems. Therefore, it is critical for us to extract and exploit any structures present in the matrix $AMA^T$ to dramatically reduce the cost of solving (29).

Next, we analyse the proof in Theorem 3.1 in details in order to find the special structure of $AMA^T$, thereby reducing the computational cost mentioned above. Let $v := \text{Prox}_\varphi(u)$. From the proof in Theorem 3.1, case 1 and case 2 (taking $t = 0$) are simple since the set $\mathcal{M}(u)$ contains a zero matrix. We can choose $M = 0$ so that
$$I + \sigma A M A^T = I.$$

The sole challenge lies in case 3. Here, we shall consider
$$\left(1 - \frac{\lambda_{2,l}}{\|v_l\|}\right) A \mathcal{P}_l^* \mathcal{P}_l \Theta A^T + \frac{\lambda_{2,l}}{\|v_l\|^3} A (\mathcal{P}_l^* v_l)(\mathcal{P}_l^* v_l)^T A^T.$$

Note that both $\mathcal{P}_l^* \mathcal{P}_l$ and $\Theta$ are diagonal matrices whose diagonal elements are either 0 or 1. Therefore, the product $\mathcal{P}_l^* \mathcal{P}_l \Theta$ enjoys the same property. Moreover, we have $\text{supp}(\text{diag}(\mathcal{P}_l^* \mathcal{P}_l)) = G_l$ and $\text{supp}(\text{diag}(\Theta)) = \text{supp}(v)$ by the definition of $\mathcal{P}_l$, (13), and (14). Therefore,
$$\text{supp}(\text{diag}(\mathcal{P}_l^* \mathcal{P}_l \Theta)) = G_l \cap \text{supp}(v) = \Xi_l,$$
where $\Xi_l$ is the index set defined by (16) that corresponds to the non-zero elements of $v$ in the $l$-th group. In other words, the diagonal matrix $\mathcal{P}_l^* \mathcal{P}_l \Theta$ is expected to contain only a few 1's in the diagonal. Consequently, the computational cost of $A \mathcal{P}_l^* \mathcal{P}_l \Theta A^T$ can be greatly reduced. Next, we observe that $\text{supp}(\mathcal{P}_l^* v_l) \subseteq \Xi_l$. Thus to compute $A(\mathcal{P}_l^* v_l)$, one just needs to consider those columns of $A$ corresponding to the index set $\Xi_l$, thereby reducing the cost of computing $A(\mathcal{P}_l^* v_l)$ and that of $A(\mathcal{P}_l^* v_l)(\mathcal{P}_l^* v_l)^T A^T$. The following notations are introduced to express these techniques clearly. Denote the index set $\Xi_> := \{l \mid \|v_l\| > \lambda_{2,l}, l = 1, 2, \ldots, g\}$, which corresponds to case 3 in Theorem 3.1. For each $l = 1, 2, \ldots, g$, let $A_l \in \mathbb{R}^{m \times |\Xi_l|}$ be the sub-matrix of $A$ with those columns in $\Xi_l$ and $s_l := (\mathcal{P}_l^* v_l)_{\Xi_l} \in \mathbb{R}^{|\Xi_l|}$ be the sub-vector of $\mathcal{P}_l^* v_l$ restricted to $\Xi_l$. Then, we deduce that

$$AMA^T = \sum_{l \in \Xi_>} \left(1 - \frac{\lambda_{2,l}}{\|v_l\|}\right) A\mathcal{P}_l^* \mathcal{P}_l \Theta A^T + \frac{\lambda_{2,l}}{\|v_l\|^3} A(\mathcal{P}_l^* v_l)(\mathcal{P}_l^* v_l)^T A^T$$

$$= \sum_{l \in \Xi_>} \left(1 - \frac{\lambda_{2,l}}{\|v_l\|}\right) A_l A_l^T + \frac{\lambda_{2,l}}{\|v_l\|^3} (A_l s_l)(A_l s_l)^T. \tag{30}$$

Therefore, the cost of computing $AMA^T$ and that of the matrix-vector product $AMA^T d$ for any $d \in \mathbb{R}^m$ are $O(m^2(r+r_2))$ and $O(m(r+r_2))$, respectively, where $r := \sum_{l \in \Xi_>} |\Xi_l| \leq |\operatorname{supp}(v)|$ and $r_2 := |\Xi_>| \leq g$. We may refer to $r$ as the overall sparsity and $r_2$ as the group sparsity. In other words, the computational cost depends on the overall sparsity $r$, the group sparsity $r_2$, and the number of observations $m$. The number $r$ is presumably much smaller than $n$ due to the fact that $v = \operatorname{Prox}_\varphi(u)$. Besides, the number of observations $m$ is usually smaller than the number of predictors $n$ in many applications. Even if $n$ happens to be extremely large (say, larger than $10^7$), one can still solve the linear system (29) efficiently via the (sparse) Cholesky factorization as long as $r$, $r_2$, and $m$ are moderate (say, less than $10^4$).

In addition, if the optimal solution is so sparse that $r + r_2 \ll m$, then the cost of solving (29) can be reduced further. In this case, the coefficient matrix can be written as follows:

$$I + \sigma AMA^T = I + DD^T,$$

where $D = [B, C] \in \mathbb{R}^{m \times (r+r_2)}$ with $B_l := \sqrt{\sigma\left(1 - \frac{\lambda_{2,l}}{\|v_l\|}\right)} A_l \in \mathbb{R}^{m \times |\Xi_l|}$, $B := [B_l]_{l \in \Xi_>} \in \mathbb{R}^{m \times r}$, $c_l := \sqrt{\sigma \frac{\lambda_{2,l}}{\|v_l\|^3}} (A_l s_l) \in \mathbb{R}^m$ and $C = [c_l]_{l \in \Xi_>} \in \mathbb{R}^{m \times r_2}$. By the Sherman-Morrison-Woodbury formula, it holds that

$$(I + \sigma AMA^T)^{-1} = (I + DD^T)^{-1} = I - D(I + D^T D)^{-1} D^T.$$

In this case, the main cost is in computing $I + D^T D$ at $O(m(r+r_2)^2)$ operations, as well as to factorize the $r+r_2$ by $r+r_2$ matrix $I + D^T D$ at the cost of $O((r+r_2)^3)$ operations.

Based on the above arguments, one can claim that the linear system (28) in each SSN iteration can be solved efficiently at low costs. In fact based on our experience gathered from the numerical experiments in the next section, the computational costs are so low that the time taken to perform indexing operations, such as obtaining the sub-matrix $A_l$ from $A$ and the sub-vector $(\mathcal{P}_l^* v_l)_{\Xi_l}$ from $\mathcal{P}_l^* v_l$ for $l \in \Xi_>$, may become noticeably higher than the time taken to compute the matrix $AMA^T$ itself. Fortunately, the group sparsity $r_2$ generally limits the number of such indexing operations needed when computing $AMA^T$.

Note that in the unlikely event that computing the Cholesky factorization of $AMA^T$ or that of $I + D^T D$ is expensive, such as when $r + r_2$ and $m$ are both large (say more than $10^4$), one can employ the preconditioned conjugate gradient (PCG) method to solve the linear system (29) efficiently through exploiting the fast computation of the matrix-vector product $AMA^T d$ for any given vector $d$.

## 5 Numerical experiments

In this section, we compare the performance of our semismooth Newton augmented Lagrangian (SSNAL) method with the alternating direction method of multipliers (ADMM) and the state-of-the-art solver SLEP[1][27] for solving the SGLasso problem. Specifically, the function "sgLeastR" in the solver SLEP is used for comparison. ADMM was first proposed in [15, 16], and the implementation will be illustrated in section 5.1.

Since the primal problem (1) is unconstrained, it is reasonable to measure the accuracy of an approximate optimal solution $(y, z, x)$ for problem (17) and problem (1) by the relative duality gap and dual infeasibility. Specifically, let

$$\text{pobj} := \frac{1}{2}\|\mathcal{A}x - b\|^2 + \lambda_1 \|x\|_1 + \lambda_2 \sum_{l=1}^{g} w_l \|x_{G_l}\| \text{ and } \text{dobj} := -\langle b, y \rangle - \frac{1}{2}\|y\|^2$$

be the primal and dual objective function values. Then the relative duality gap and the relative dual infeasibility are defined by

$$\eta_G := \frac{|\text{pobj} - \text{dobj}|}{1 + |\text{pobj}| + |\text{dobj}|}, \quad \eta_D := \frac{\|\mathcal{A}^* y + z\|}{1 + \|z\|}.$$

For given error tolerances $\varepsilon_D > 0$ and $\varepsilon_G > 0$, our algorithm SSNAL will be terminated if

$$\eta_D < \varepsilon_D \quad \text{and} \quad \eta_G < \varepsilon_G, \tag{31}$$

while the ADMM will be terminated if the above conditions hold or the maximum number of 10,000 iterations is reached. By contrast, since SLEP does not produce the dual sequences $\{(y^k, z^k)\}$, the relative dual infeasibility cannot be used as a stopping criterion for SLEP. Therefore, We terminate SLEP if the relative difference of the optimal objective values between SLEP and SSNAL is less than $\varepsilon_G$, i.e.,

$$\eta_P := \frac{\text{obj}_S - \text{pobj}}{1 + |\text{obj}_S| + |\text{pobj}|} < \varepsilon_G,$$

or the maximum number of 10,000 iterations is reached. Here $\text{obj}_S$ denotes the objective value obtained by SLEP. Note that the parameters for SLEP are set to their default values unless otherwise specified.

---

[1]http://www.public.asu.edu/~jye02/Software/SLEP

In numerical our experiments, we choose $\varepsilon_D = \varepsilon_G = 10^{-6}$ unless otherwise specified. That is, the condition (31) for SSNAL becomes

$$\eta_S := \max\{\eta_G, \eta_D\} < 10^{-6}.$$

Similarly, the stopping condition for ADMM becomes

$$\eta_A := \max\{\eta_G, \eta_D\} < 10^{-6}.$$

In addition, we adopt the following weights: $w_l = \sqrt{|G_l|}$, $\forall l = 1, 2, \ldots, g$ for the model (1). In the following tables, "S" stands for SSNAL; "P" for SLEP; "A" for ADMM; "nnz" denotes the number of non-zero entries in the solution $x$ obtained by SSNAL using the following estimation:

$$\text{nnz} := \min\{k \mid \sum_{i=1}^{k} |\hat{x}_i| \geq 0.999\|x\|_1\},$$

where $\hat{x}$ is obtained via sorting $x$ by magnitude in a descending order. We display the number of outer ALM iterations (in Algorithm 1) and the total number of inner SSN iterations (in Algorithm 2) of SSNAL in the format of "outer iteration (inner iteration)" under the iteration column. The computation time is in the format of "hours:minutes:seconds", and "00" in the time column means that the elapsed time is less than 0.5 second.

All our numerical results are obtained by running MATLAB (version 9.0) on a windows workstation (24-core, Intel Xeon E5-2680 @ 2.50GHz, 128 Gigabytes of RAM).

## 5.1 Dual based ADMM

In this section, we study the implementation of the (inexact) semi-proximal alternating direction method of multipliers (sPADMM), which is an extension of the classic ADMM [15, 16]. This method is one of the most natural methods for solving (17) due to its separable structure. Generally, the framework of the sPADMM consists of the following iterations:

$$\begin{cases} y^{k+1} \approx \arg\min_y \mathcal{L}_\sigma(y, z^k; x^k) + \frac{1}{2}\|y - y^k\|^2_{\mathcal{S}_1}, \\ z^{k+1} \approx \arg\min_z \mathcal{L}_\sigma(y^{k+1}, z; x^k) + \frac{1}{2}\|z - z^k\|^2_{\mathcal{S}_2}, \\ x^{k+1} = x^k - \tau\sigma(\mathcal{A}^*y^{k+1} + z^{k+1}), \end{cases} \quad (32)$$

where $\tau \in (0, (1+\sqrt{5})/2)$, $\mathcal{S}_1$ and $\mathcal{S}_2$ are self-adjoint positive semidefinite linear operators, and $\mathcal{L}_\sigma$ is the augmented Lagrangian function defined in (18). The sPADMM is convergent under some mild conditions, and we refer the reader to [12, 4] for the convergence results. However, due to the lack of error bound conditions for the KKT system (2), the linear convergence rate of the sPADMM cannot be established from existing results.

In each iteration of (32), the first step is to minimize a function of $y$. In particular, $y^{k+1}$ can be obtained by solving the following $m \times m$ linear system of equations:

$$(\sigma^{-1}I + \mathcal{A}\mathcal{A}^* + \mathcal{S}_1)y^{k+1} = -\sigma^{-1}b - \mathcal{A}(z^k - \sigma^{-1}x^k) + \mathcal{S}_1 y^k.$$

As the dimension $m$ is a moderate number in many statistical applications. Thus, in our implementation, equation (5.1) was solved via the Cholesky factorization, and the proximal term $\mathcal{S}_1$ was taken to be the zero matrix. In the event that computing the Cholesky factorization of $\sigma^{-1}I + \mathcal{A}\mathcal{A}^*$ is expensive, one can choose $\mathcal{S}_1$ judiciously to make the coefficient matrix to be a positive definite diagonal matrix plus a low-rank matrix for which one can invert it efficiently via the Sherman-Morrison-Woodbury formula. We refer the reader to [4, section 7.1] for the details on how to choose $\mathcal{S}_1$ appropriately.

The second step in (32) is to minimize a function of $z$. For the SGLasso problem, one would simply choose $\mathcal{S}_2 = 0$. In this case, by the Moreau identity (4), $z^{k+1}$ is updated by the following scheme:

$$z^{k+1} = \sigma^{-1}x^k - \mathcal{A}^*y^{k+1} - \text{Prox}_p(\sigma^{-1}x^k - \mathcal{A}^*y^{k+1}),$$

where $\text{Prox}_p$ is computable by Proposition 2.1. In summary, two subproblems of (32) are solvable and consequently the framework (32) is easily implementable. Moreover, in order to improve the convergence speed numerically, we set the step-length $\tau$ in (32) to be 1.618 and tune the parameter $\sigma$ according to the progress between primal feasibility and dual feasibility in the implementation.

## 5.2 Synthetic data

This section presents the tests of the three algorithms SSNAL, ADMM, and SLEP on various synthetic data constructed in the same way as in [38]. The data matrix $A$ is generated randomly as an $m \times n$ matrix of normally distributed random numbers, and the number of groups $g$ is chosen manually to be 100, 1000, and 10000. Then we partition $\{1, 2, \ldots, n\}$ into $g$ groups such that the indices of components in each group are adjacent, for example, $G_1 = \{1, 2, \ldots, 25\}$, $G_2 = \{26, 27, \ldots, 53\}$, etc. The group sizes $\{|G_i|, i = 1, 2, \ldots, g\}$ are determined randomly such that each $|G_i|$ is expected to be around the mean value of $\frac{n}{g}$. Subsequently, the response vector $b$ is constructed as

$$b = Ax + \epsilon,$$

where $\epsilon$ is normally distributed random noise, $x_{G_l} = (1, 2, \ldots, 10, 0, \ldots, 0)^T$ for $l = 1, 2, \ldots, 10$, and $x_{G_l} = 0$ for all other groups. That is, the first 10 groups are the non-trivial groups, and the true number of non-zero elements of the underlying solution $x$ is 100. The regularization parameters $\lambda_1 = \lambda_2$ are chosen to make the number of non-zero elements of the resulting solution close to the true number of 100.

Table 1 compares the numerical results of the three algorithms SSNAL, ADMM, and SLEP tested on different synthetic data. As can be seen from the table, the computational time of SSNAL is less than that of ADMM and SLEP for most cases. The overall advantage of computational time suggests that our algorithm SSNAL is efficient for solving the SGLasso problem with randomly generated data. Moreover, we observe from the table that ADMM is inefficient in solving the SGLasso problem with randomly generated large-scale data. A possible reason is that the first order method ADMM requires a large number of iterations to solve the problem to the required accuracy of $10^{-6}$. The table also shows that our algorithm SSNAL can significantly outperform SLEP on problems with a large number of groups. In particular, SSNAL is more than 5 times faster than SLEP for the high dimensional instance with problem size $(m,n) = (1e4, 1e6)$ and group number $g = 10000$. For this instance, the number of non-zero entries in the solution $x$ is small, and we have highly conducive second order sparsity which we can fully exploit in the numerical computations outlined in section 4.3.

Table 1: The performances of SSNAL, ADMM, and SLEP on synthetic data. Regularization parameters are set as follows: $\lambda_1 = \lambda_2$.

| size $(m,n)$ | $g$ | $\lambda_1$ | nnz | iteration S\|A\|P | time S\|A\|P |
|---|---|---|---|---|---|
| (1e3,1e5) | 100 | 1338 | 166 | 1(3) \| 1246 \| 1 | 00 \| 01:23 \| 00 |
|  | 1000 | 1736 | 154 | 2(13) \| 1247 \| 26 | 01 \| 01:25 \| 01 |
|  | 10000 | 983 | 84 | 4(27) \| 1185 \| 239 | 02 \| 01:21 \| 13 |
| (1e4,1e6) | 100 | 3775 | 43 | 1(3) \| 2228 \| 17 | 13 \| 03:01:49 \| 59 |
|  | 1000 | 7229 | 167 | 1(3) \| 2232 \| 1 | 11 \| 03:05:03 \| 08 |
|  | 10000 | 4000 | 109 | 3(28) \| 2104 \| 148 | 01:38 \| 03:06:31 \| 08:37 |

## 5.3 UCI data sets with random groups

This section presents the performances of the three algorithms SSNAL, ADMM, and SLEP on large-scale UCI data sets [26] $(\mathcal{A}, b)$ that are originally obtained from the LIBSVM data sets [3]. In our numerical experiments, we follow [24] and apply the method in [18] to expand the original features of the data sets *bodyfat*, *pyrim*, and *triazines* using polynomial basis functions. For example, a polynomial basis function of order 7 is used to expand the features of the data set *bodyfat*, and then the expanded data set is named as *bodyfat7*. This naming convention is also used for *pyrim5*, *triazines4*, and *housing 7*. As noted in [24, Table 1], these data sets are quite different in terms of the problem dimension and the largest eigenvalue of $\mathcal{A}\mathcal{A}^*$. For example, for a relatively high-dimensional instance *log1p.E2006.train*, the dimension of $\mathcal{A}$ is $16087 \times 4272227$ and the largest eigenvalue of $\mathcal{A}\mathcal{A}^*$ is $5.86 \times 10^7$.

Next, we describe how the groups in each problem are specified. By reordering the components of the variable $x$ if necessary, without loss of generality, we assume that the vector $x$ can be partitioned into $g$ groups where the indices of components in each group are adjacent. The group sizes $\{|G_l|, l = 1, 2, \ldots, g\}$ are determined randomly such that each $|G_l|$ is around the mean value of $\frac{n}{g}$. In the experiment, the average group size is about 300.

We tested the SGLasso problems with two different sets of regularization parameters which are chosen manually:

$$\text{(S1)} \quad \lambda_1 = \lambda_2 = \gamma\|\mathcal{A}^*b\|_\infty;$$

$$\text{(S2)} \quad \lambda_1 = 0.5\gamma\|\mathcal{A}^*b\|_\infty, \ \lambda_2 = 9.5\gamma\|\mathcal{A}^*b\|_\infty.$$

The parameter $\gamma$ is chosen to produce a reasonable number of non-zero elements in the resulting solution $x$. Three values of $\gamma$ are used for each UCI data set in our experiments.

Table 2 presents the comparison results of the three algorithms Ssnal, ADMM, and SLEP on 8 selected UCI data sets with regularization parameters specified as in (S1). As shown in the table, Ssnal has succeeded in solving all instances within 1 minutes, while SLEP failed to solve 10 cases. Although ADMM has also succeeded in solving all instances, its running time for each case is much longer than that of Ssnal. In majority of the cases, Ssnal outperformed the first order methods ADMM and SLEP by a large margin. For example, for the instance *E2006.train* with $\gamma = 1e$-7, Ssnal solved it to the desired accuracy in 3 seconds, ADMM took more than 8 minutes, while SLEP failed to solve it within 10000 steps. The numerical results show convincingly that our algorithm Ssnal can solve SGLasso problems highly efficiently and robustly. Again, the superior performance of our Ssnal algorithm can be attributed to our ability to extract and exploit the second order sparsity structure (in the SGLasso problem) within the SSN method to solve each ALM subproblem very efficiently.

Table 3 is the same as Table 2 but for the regularization parameters specified as in (S2). This table also shows that the computational time of Ssnal is far less than that of ADMM and SLEP for almost all cases. Furthermore, for more difficult cases, such as those with large problem dimension $(m,n)$ and large number of non-zero entries (nnz), the superiority of Ssnal is even more striking compared to ADMM and SLEP. The results again demonstrate that our algorithm Ssnal is highly efficient for solving SGLasso problems.

Figure 1 presents the performance profiles of Ssnal, ADMM, and SLEP for all 48 tested problems, which are presented in Table 2 and Table 3. The meaning of the performance profiles is given as follows: a point $(x, y)$ is on the performance curve of a particular method if and only if this method can solve up to desired accuracy $(100y)\%$ of all the tested instances within at most $x$ times of the fastest method for each instance. As can be seen, Ssnal outperforms ADMM and SLEP by a large margin for all tested UCI data sets with randomly generated groups. In particular, focusing on $y = 40\%$, we can see from Figure 1 that Ssnal is around 30 times faster compared to ADMM and SLEP for over 60% of the tested instances.

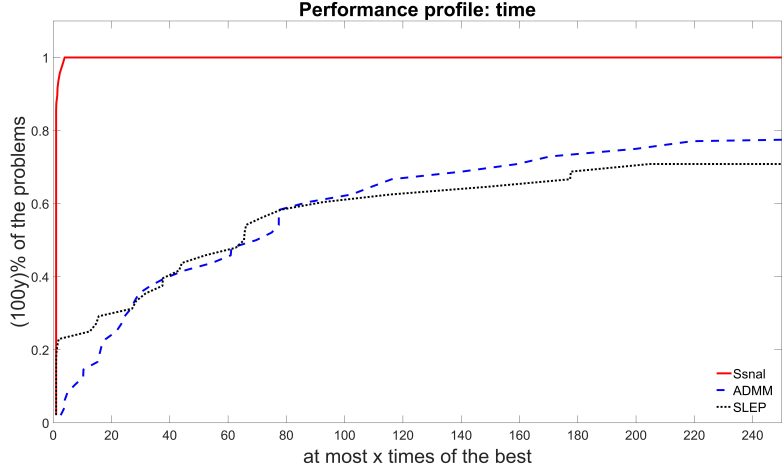

Figure 1: Performance profiles of Ssnal, ADMM, and SLEP on UCI data sets with randomly generated groups.

Table 2: The performances of Ssnal, ADMM, and SLEP on 8 selected UCI data sets with randomly generated groups. The regularization parameters are specified as in $(S_1)$.

| problem name (m,n); g | $\gamma$ | nnz | iteration S\|A\|P | time S\|A\|P | error $\eta_S\|\eta_A\|\eta_P$ |
|---|---|---|---|---|---|
| E2006.train     | 1e-05 | 1    | 3(7) \| 44 \| 4        | 01 \| 07:50 \| 00    | 1.4e-09 \| 7.4e-07 \| 6.2e-09 |
| (16087,150360)  | 1e-06 | 1    | 4(8) \| 23 \| 11       | 01 \| 07:03 \| 01    | 7.9e-11 \| 7.8e-07 \| 9.9e-07 |
| 501             | 1e-07 | 46   | 20(40) \| 70 \| 10000  | 03 \| 08:30 \| 08:59 | 1.2e-07 \| 9.3e-07 \| 2.5e-04 |
| E2006.test      | 1e-05 | 1    | 4(8) \| 36 \| 4        | 00 \| 18 \| 00       | 1.0e-10 \| 8.6e-07 \| 4.9e-13 |
| (3308,150358)   | 1e-06 | 1    | 4(9) \| 29 \| 21       | 00 \| 18 \| 00       | 1.2e-08 \| 8.3e-07 \| 9.7e-07 |
| 501             | 1e-07 | 150  | 23(48) \| 205 \| 10000 | 02 \| 28 \| 02:17    | 8.0e-07 \| 1.0e-06 \| 3.1e-03 |
| log1p.E2006.train | 1e-03 | 2   | 3(10) \| 2360 \| 882   | 09 \| 01:01:06 \| 04:45 | 4.9e-09 \| 1.0e-06 \| 9.7e-07 |
| (16087,4272227) | 1e-04 | 3    | 3(10) \| 2490 \| 5260  | 08 \| 01:02:49 \| 28:21 | 7.0e-07 \| 1.0e-06 \| 9.8e-07 |
| 14241           | 1e-05 | 1005 | 5(22) \| 876 \| 7553   | 37 \| 33:57 \| 40:48 | 3.3e-07 \| 9.9e-07 \| 1.0e-06 |
| log1p.E2006.test | 1e-03 | 4   | 4(13) \| 1549 \| 2107  | 08 \| 10:22 \| 05:12 | 3.5e-07 \| 9.9e-07 \| 9.8e-07 |
| (3308,4272226)  | 1e-04 | 5    | 4(12) \| 1749 \| 2693  | 06 \| 11:39 \| 06:41 | 6.2e-07 \| 9.9e-07 \| 9.9e-07 |
| 14241           | 1e-05 | 5009 | 7(34) \| 464 \| 10000  | 45 \| 03:42 \| 25:04 | 2.1e-07 \| 9.7e-07 \| 7.5e-06 |
| bodyfat7        | 1e-04 | 7    | 11(29) \| 878 \| 3246  | 01 \| 28 \| 54       | 7.7e-07 \| 9.9e-07 \| 9.6e-07 |
| (252,116280)    | 1e-05 | 13   | 15(37) \| 918 \| 10000 | 03 \| 29 \| 02:48    | 7.9e-07 \| 9.9e-07 \| 2.1e-05 |
| 388             | 1e-06 | 237  | 21(58) \| 917 \| 10000 | 07 \| 29 \| 02:51    | 6.9e-07 \| 1.0e-06 \| 2.4e-04 |
| pyrim5          | 1e-02 | 279  | 8(32) \| 3074 \| 4736  | 02 \| 02:04 \| 01:55 | 1.3e-07 \| 1.0e-06 \| 1.0e-06 |
| (74,201376)     | 1e-03 | 606  | 11(39) \| 2003 \| 10000 | 02 \| 01:21 \| 04:05 | 3.9e-07 \| 1.0e-06 \| 8.9e-06 |
| 671             | 1e-04 | 937  | 17(50) \| 1969 \| 10000 | 05 \| 01:23 \| 04:05 | 6.1e-07 \| 1.0e-06 \| 1.2e-04 |
| triazines4      | 1e-02 | 406  | 8(35) \| 5038 \| 9374  | 09 \| 24:37 \| 25:38 | 8.5e-08 \| 1.0e-06 \| 1.0e-06 |
| (186,635376)    | 1e-03 | 1396 | 9(43) \| 4020 \| 10000 | 17 \| 19:43 \| 27:25 | 4.2e-07 \| 1.0e-06 \| 5.8e-04 |
| 2118            | 1e-04 | 3574 | 16(58) \| 4287 \| 10000 | 55 \| 22:33 \| 27:23 | 6.6e-07 \| 1.0e-06 \| 5.5e-03 |
| housing7        | 1e-02 | 220  | 7(31) \| 813 \| 3366   | 01 \| 25 \| 59       | 2.5e-07 \| 9.9e-07 \| 9.9e-07 |
| (506,77520)     | 1e-03 | 817  | 9(37) \| 816 \| 5199   | 02 \| 25 \| 01:32    | 8.9e-08 \| 9.9e-07 \| 1.0e-06 |
| 258             | 1e-04 | 2134 | 14(47) \| 618 \| 10000 | 07 \| 19 \| 02:56    | 7.1e-07 \| 1.0e-06 \| 3.7e-06 |

Table 3: The performances of SSNAL, ADMM, and SLEP on 8 selected UCI data sets with randomly generated groups. The regularization parameters are specified as in $(S_2)$.

| problem name (m,n);g | $\gamma$ | nnz | iteration S\|A\|P | time S\|A\|P | error $\eta_S\|\eta_A\|\eta_P$ |
|---|---|---|---|---|---|
| E2006.train (16087,150360) 501 | 1e-05 | 1 | 3(7) \| 73 \| 1 | 01 \| 08:40 \| 00 | 1.6e-10 \| 9.7e-07 \| 3.0e-08 |
|  | 1e-06 | 1 | 3(7) \| 44 \| 4 | 01 \| 07:49 \| 00 | 1.5e-09 \| 7.9e-07 \| 7.3e-07 |
|  | 1e-07 | 16 | 8(16) \| 31 \| 13 | 01 \| 07:33 \| 01 | 3.3e-07 \| 9.9e-07 \| 9.9e-07 |
| E2006.test (3308,150358) 501 | 1e-05 | 1 | 3(7) \| 52 \| 16 | 00 \| 20 \| 00 | 3.7e-09 \| 7.1e-07 \| 3.8e-08 |
|  | 1e-06 | 1 | 4(8) \| 32 \| 16 | 00 \| 18 \| 00 | 4.9e-10 \| 9.4e-07 \| 4.9e-08 |
|  | 1e-07 | 15 | 9(19) \| 29 \| 25 | 01 \| 18 \| 00 | 4.8e-07 \| 6.9e-07 \| 8.7e-07 |
| log1p.E2006.train (16087,4272227) 14241 | 1e-03 | 1 | 2(6) \| 2664 \| 202 | 05 \| 01:06:05 \| 01:06 | 5.1e-08 \| 9.9e-07 \| 2.0e-07 |
|  | 1e-04 | 7 | 3(11) \| 2379 \| 1286 | 10 \| 01:00:42 \| 06:53 | 4.2e-09 \| 9.9e-07 \| 9.6e-07 |
|  | 1e-05 | 32 | 3(11) \| 2474 \| 4375 | 09 \| 01:02:34 \| 23:23 | 6.8e-07 \| 1.0e-06 \| 9.9e-07 |
| log1p.E2006.test (3308,4272226) 14241 | 1e-03 | 2 | 2(7) \| 1961 \| 379 | 04 \| 12:55 \| 57 | 5.5e-07 \| 9.9e-07 \| 9.0e-07 |
|  | 1e-04 | 10 | 4(13) \| 1567 \| 1459 | 08 \| 10:30 \| 03:40 | 3.9e-07 \| 9.9e-07 \| 9.6e-07 |
|  | 1e-05 | 95 | 5(15) \| 1749 \| 4800 | 08 \| 11:37 \| 12:18 | 5.6e-08 \| 9.9e-07 \| 9.9e-07 |
| bodyfat7 (252,116280) 388 | 1e-04 | 111 | 9(20) \| 363 \| 1711 | 00 \| 12 \| 33 | 1.7e-08 \| 9.8e-07 \| 1.0e-06 |
|  | 1e-05 | 208 | 13(30) \| 438 \| 8460 | 01 \| 14 \| 02:42 | 2.5e-07 \| 9.6e-07 \| 1.0e-06 |
|  | 1e-06 | 264 | 17(37) \| 555 \| 10000 | 05 \| 18 \| 03:10 | 7.9e-07 \| 9.9e-07 \| 2.3e-06 |
| pyrim5 (74,201376) 671 | 1e-02 | 230 | 4(17) \| 1258 \| 1489 | 01 \| 51 \| 37 | 3.7e-07 \| 4.5e-07 \| 9.7e-07 |
|  | 1e-03 | 626 | 8(34) \| 1413 \| 6038 | 02 \| 57 \| 02:32 | 9.0e-07 \| 1.0e-06 \| 1.0e-06 |
|  | 1e-04 | 1178 | 13(43) \| 1684 \| 10000 | 04 \| 01:10 \| 04:16 | 1.1e-07 \| 1.0e-06 \| 4.8e-05 |
| triazines4 (186,635376) 2118 | 1e-02 | 577 | 6(27) \| 10000 \| 4422 | 06 \| 48:29 \| 12:01 | 1.1e-07 \| 2.7e-06 \| 1.0e-06 |
|  | 1e-03 | 1171 | 8(36) \| 4875 \| 10000 | 10 \| 23:49 \| 27:19 | 5.3e-07 \| 1.0e-06 \| 3.2e-06 |
|  | 1e-04 | 4346 | 11(48) \| 3343 \| 10000 | 28 \| 17:51 \| 28:42 | 9.1e-07 \| 1.0e-06 \| 2.7e-04 |
| housing7 (506,77520) 258 | 1e-02 | 206 | 3(11) \| 1097 \| 29 | 00 \| 34 \| 01 | 9.0e-09 \| 9.7e-07 \| 2.9e-07 |
|  | 1e-03 | 839 | 8(30) \| 754 \| 936 | 01 \| 23 \| 17 | 1.3e-07 \| 1.0e-06 \| 9.8e-07 |
|  | 1e-04 | 1689 | 10(36) \| 837 \| 5510 | 03 \| 26 \| 01:38 | 1.0e-07 \| 1.0e-06 \| 1.0e-06 |

## 5.4 UCI datesets with simulated groups

This section also makes uses of the UCI data sets mentioned in section 5.3. Instead of specifying the groups randomly, we attempt to generate more meaningful groups in the following manner. Firstly, the classical Lasso (model (1) with $\lambda_2 = 0$) is solved with the accuracy of $10^{-4}$ to obtain a sparse solution $x$, and the computed solution $x$ is sorted in a descending order. Then, the first $|G_1|$ largest variables are allocated to group 1, and the next $|G_2|$ variables are allocated to group 2, etc. Since this group membership is determined by the magnitude of each variable of the computed solution from the classical Lasso, we believe that this kind of group structure is more natural than that constructed randomly in the previous section. Besides, the group sizes $\{|G_l|, l = 1, 2, \ldots, g\}$ are determined randomly such that each $|G_l|$ is around the mean value of $\frac{n}{g}$. Compared to the last section, a different value 30 is taken as the average group size for the diversity of experiments.

To generate the solution from the classical Lasso to decide on the group membership mentioned above, we take the medium value of $\gamma$ in Table 4, e.g., $\gamma = 1e$-6 for the instance *E2006.train*. And the regularization parameters for the classical Lasso are set as follow: $\lambda_1 = \gamma \|\mathcal{A}^* b\|_\infty$, $\lambda_2 = 0$. For the SGLasso problem, the regularization parameters follow three different strategies: (S1) and (S2) given in the previous section, and

$$(S3) \quad \lambda_1 = \gamma \|\mathcal{A}^* b\|_\infty, \ \lambda_2 = \sqrt{\lambda_1} \text{ if } \lambda_1 > 1 \text{ and } \lambda_2 = \lambda_1^2 \text{ if } \lambda_1 \leq 1.$$

The comparison results with parameter sets (S1), (S2), and (S3) are presented in Table 4, Table 5, and Table 6, respectively. As shown in these three tables, S$\textsc{snal}$ has succeeded in solving all the 72 instances highly efficiently, while ADMM failed in 5 instances, and SLEP failed in 58 instances. Moreover, for those failed instances, we observe from the tables that SLEP terminated when the errors are still relatively large, which is $10^{-2}$ for most cases. The results may suggest that using only first order information is not enough for computing high accuracy solution, while second order information can contribute to the fast convergence and high computational efficiency of a well designed second order SSN method. For the vast majority of the instances, the computational time of S$\textsc{snal}$ is far less than that of ADMM and SLEP. Again, the results have demonstrated convincingly that our algorithm S$\textsc{snal}$ is capable of solving large-scale SGLasso problems to high accuracy very efficiently and robustly.

Figure 2 presents the performance profiles of S$\textsc{snal}$, ADMM, and SLEP for all 72 tested problems, which are presented in Table 4, Table 5, and Table 6. From the figure, we find that S$\textsc{snal}$ not only solves all the tested instances to the desired accuracy, but also outperforms ADMM and SLEP by an obvious margin for these tested UCI data sets with simulated groups. Within 250 times of the running time of S$\textsc{snal}$, ADMM can only solve approximately 80% of all the tested instances, while SLEP can only solve 20% of all the tested instances. We can safely claim that our algorithm S$\textsc{snal}$ can solve large-scale SGLasso problems to high accuracy very efficiently and robustly.

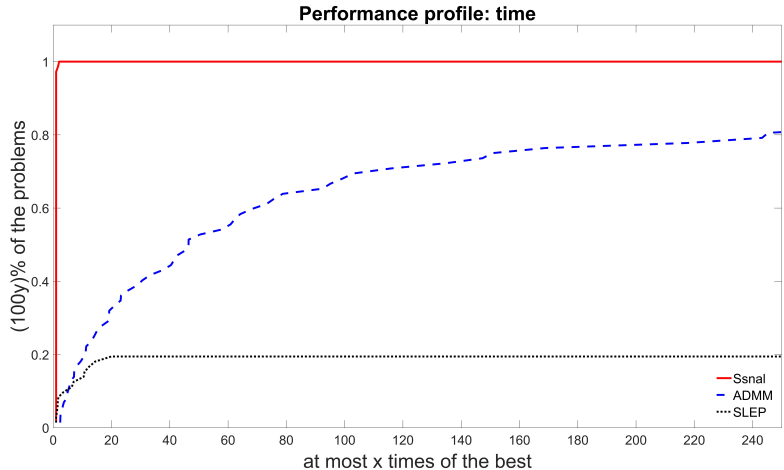

Figure 2: Performance profiles of S$\textsc{snal}$, ADMM, and SLEP on UCI data sets with simulated groups.

Table 4: The performances of SSNAL, ADMM, and SLEP on 8 selected UCI data sets with simulated groups. The regularization parameters are specified as in $(S_1)$.

| problem name (m,n);g | $\gamma$ | nnz | iteration S\|A\|P | time S\|A\|P | error $\eta_S\|\eta_A\|\eta_P$ |
|---|---|---|---|---|---|
| E2006.train | 1e-05 | 1 | 4(9) \| 36 \| 16 | 01 \| 07:30 \| 01 | 6.4e-10 \| 9.0e-07 \| 4.1e-08 |
| (16087,150360) | 1e-06 | 34 | 14(28) \| 33 \| 10000 | 02 \| 07:14 \| 08:56 | 9.2e-07 \| 6.9e-07 \| 3.7e-05 |
| 5012 | 1e-07 | 210 | 25(52) \| 110 \| 10000 | 08 \| 09:17 \| 08:58 | 5.8e-08 \| 7.7e-07 \| 1.2e-02 |
| E2006.test | 1e-05 | 1 | 4(9) \| 31 \| 9 | 00 \| 17 \| 00 | 5.4e-09 \| 9.7e-07 \| 7.7e-07 |
| (3308,150358) | 1e-06 | 38 | 20(41) \| 70 \| 10000 | 01 \| 20 \| 02:35 | 3.1e-07 \| 8.2e-07 \| 1.3e-03 |
| 5012 | 1e-07 | 275 | 27(56) \| 324 \| 10000 | 02 \| 33 \| 02:38 | 1.9e-07 \| 9.5e-07 \| 9.8e-02 |
| log1p.E2006.train | 1e-03 | 15 | 3(13) \| 2249 \| 290 | 09 \| 57:09 \| 01:37 | 5.5e-07 \| 9.9e-07 \| -1.5e-04 |
| (16087,4272227) | 1e-04 | 253 | 4(25) \| 1194 \| 10000 | 29 \| 38:34 \| 54:46 | 1.4e-07 \| 1.0e-06 \| 1.7e-02 |
| 142408 | 1e-05 | 7821 | 6(32) \| 391 \| 10000 | 03:27 \| 24:39 \| 54:48 | 2.4e-07 \| 9.8e-07 \| 1.5e-02 |
| log1p.E2006.test | 1e-03 | 13 | 4(14) \| 1572 \| 284 | 07 \| 10:16 \| 47 | 4.6e-07 \| 9.9e-07 \| -2.4e-04 |
| (3308,4272226) | 1e-04 | 546 | 5(21) \| 627 \| 10000 | 14 \| 04:33 \| 26:56 | 2.1e-07 \| 9.9e-07 \| 4.0e-02 |
| 142408 | 1e-05 | 4874 | 8(33) \| 300 \| 10000 | 43 \| 02:37 \| 27:04 | 6.2e-07 \| 9.8e-07 \| 1.1e-01 |
| bodyfat7 | 1e-04 | 11 | 12(32) \| 930 \| 10000 | 01 \| 30 \| 02:54 | 8.1e-07 \| 9.9e-07 \| 6.7e-04 |
| (252,116280) | 1e-05 | 26 | 19(53) \| 2394 \| 10000 | 03 \| 01:16 \| 02:56 | 4.1e-07 \| 1.0e-06 \| 2.2e-04 |
| 3876 | 1e-06 | 166 | 23(75) \| 1201 \| 10000 | 08 \| 38 \| 02:58 | 2.5e-08 \| 9.9e-07 \| 2.0e-04 |
| pyrim5 | 1e-02 | 98 | 7(27) \| 1243 \| 10000 | 02 \| 51 \| 03:39 | 3.6e-07 \| 9.9e-07 \| 9.4e-02 |
| (74,201376) | 1e-03 | 201 | 12(43) \| 2080 \| 10000 | 02 \| 01:27 \| 03:39 | 2.1e-07 \| 1.0e-06 \| 4.2e-02 |
| 6713 | 1e-04 | 644 | 18(66) \| 2351 \| 10000 | 06 \| 01:43 \| 03:46 | 3.4e-07 \| 1.0e-06 \| 1.0e-02 |
| triazines4 | 1e-02 | 261 | 10(42) \| 8439 \| 10000 | 11 \| 44:25 \| 26:35 | 4.1e-08 \| 9.6e-07 \| 6.8e-02 |
| (186,635376) | 1e-03 | 737 | 15(62) \| 10000 \| 10000 | 18 \| 50:16 \| 26:36 | 3.4e-08 \| 1.2e-05 \| 6.5e-02 |
| 21179 | 1e-04 | 1510 | 20(79) \| 9466 \| 10000 | 36 \| 01:20:53 \| 39:04 | 4.3e-08 \| 1.0e-06 \| 5.7e-02 |
| housing7 | 1e-02 | 91 | 6(25) \| 862 \| 10000 | 01 \| 30 \| 03:00 | 3.6e-08 \| 9.9e-07 \| 2.5e-02 |
| (506,77520) | 1e-03 | 150 | 9(34) \| 596 \| 10000 | 02 \| 21 \| 03:01 | 8.0e-07 \| 9.9e-07 \| 9.4e-02 |
| 2584 | 1e-04 | 807 | 15(49) \| 638 \| 10000 | 09 \| 23 \| 03:01 | 6.0e-07 \| 1.0e-06 \| 3.8e-02 |

Table 5: The performances of SSNAL, ADMM, and SLEP on 8 selected UCI data sets with simulated groups. The regularization parameters are specified as in $(S_2)$.

| problem name (m,n);g | $\gamma$ | nnz | iteration S\|A\|P | time S\|A\|P | error $\eta_S\|\eta_A\|\eta_P$ |
|---|---|---|---|---|---|
| E2006.train | 1e-05 | 1 | 3(7) \| 54 \| 16 | 01 \| 08:02 \| 01 | 6.4e-10 \| 7.1e-07 \| 3.9e-08 |
| (16087,150360) | 1e-06 | 11 | 8(17) \| 34 \| 10000 | 01 \| 07:28 \| 08:46 | 3.0e-07 \| 8.6e-07 \| 1.6e-05 |
| 5012 | 1e-07 | 40 | 21(47) \| 68 \| 10000 | 03 \| 08:15 \| 08:54 | 8.7e-08 \| 9.9e-07 \| 3.3e-03 |
| E2006.test | 1e-05 | 2 | 5(11) \| 42 \| 10000 | 00 \| 18 \| 02:08 | 2.9e-07 \| 7.0e-07 \| 1.5e-06 |
| (3308,150358) | 1e-06 | 22 | 9(19) \| 29 \| 10000 | 01 \| 17 \| 02:26 | 3.3e-07 \| 7.3e-07 \| 6.6e-05 |
| 5012 | 1e-07 | 66 | 25(51) \| 138 \| 10000 | 01 \| 24 \| 02:30 | 5.1e-07 \| 8.4e-07 \| 1.5e-02 |
| log1p.E2006.train | 1e-03 | 11 | 2(12) \| 2399 \| 92 | 09 \| 01:00:29 \| 30 | 6.5e-08 \| 1.0e-06 \| -2.9e-04 |
| (16087,4272227) | 1e-04 | 39 | 3(16) \| 2088 \| 578 | 13 \| 55:02 \| 03:14 | 4.5e-07 \| 1.0e-06 \| -1.8e-05 |
| 142408 | 1e-05 | 597 | 4(22) \| 862 \| 10000 | 29 \| 33:09 \| 55:12 | 2.7e-07 \| 9.9e-07 \| 3.2e-02 |
| log1p.E2006.test | 1e-03 | 7 | 2(12) \| 1567 \| 60 | 07 \| 10:20 \| 10 | 8.6e-07 \| 1.0e-06 \| -4.4e-04 |
| (3308,4272226) | 1e-04 | 47 | 4(17) \| 1260 \| 327 | 08 \| 08:28 \| 53 | 1.3e-07 \| 1.0e-06 \| -1.3e-04 |
| 142408 | 1e-05 | 1079 | 5(23) \| 467 \| 10000 | 17 \| 03:37 \| 27:49 | 9.8e-07 \| 9.9e-07 \| 1.2e-01 |
| bodyfat7 | 1e-04 | 26 | 10(24) \| 748 \| 10000 | 01 \| 24 \| 03:23 | 3.7e-07 \| 1.0e-06 \| 2.1e-02 |
| (252,116280) | 1e-05 | 43 | 15(37) \| 1266 \| 10000 | 01 \| 41 \| 03:22 | 6.1e-07 \| 1.0e-06 \| 2.4e-03 |
| 3876 | 1e-06 | 52 | 20(53) \| 1188 \| 10000 | 04 \| 38 \| 03:25 | 2.4e-07 \| 1.0e-06 \| 3.9e-04 |
| pyrim5 | 1e-02 | 42 | 6(19) \| 1672 \| 10000 | 01 \| 01:05 \| 04:02 | 6.4e-08 \| 1.0e-06 \| 1.1e-01 |
| (74,201376) | 1e-03 | 136 | 8(32) \| 1518 \| 10000 | 01 \| 59 \| 04:24 | 1.5e-07 \| 9.9e-07 \| 1.2e-01 |
| 6713 | 1e-04 | 342 | 13(50) \| 1879 \| 10000 | 04 \| 01:49 \| 04:27 | 1.6e-07 \| 1.0e-06 \| 3.7e-02 |
| triazines4 | 1e-02 | 40 | 8(20) \| 6085 \| 10000 | 04 \| 30:37 \| 26:40 | 1.6e-08 \| 9.0e-07 \| 1.1e-01 |
| (186,635376) | 1e-03 | 544 | 10(43) \| 6473 \| 10000 | 11 \| 32:06 \| 26:29 | 7.4e-08 \| 9.7e-07 \| 7.0e-02 |
| 21179 | 1e-04 | 964 | 17(63) \| 10000 \| 10000 | 18 \| 49:47 \| 26:52 | 4.1e-08 \| 2.2e-06 \| 8.2e-02 |
| housing7 | 1e-02 | 51 | 4(15) \| 1242 \| 10000 | 00 \| 38 \| 03:36 | 5.4e-08 \| 9.8e-07 \| 5.1e-02 |
| (506,77520) | 1e-03 | 153 | 7(26) \| 853 \| 10000 | 01 \| 26 \| 03:36 | 1.2e-07 \| 1.0e-06 \| 5.0e-02 |
| 2584 | 1e-04 | 175 | 10(34) \| 577 \| 10000 | 02 \| 18 \| 03:34 | 1.2e-07 \| 9.8e-07 \| 1.3e-01 |

Table 6: The performances of SSNAL, ADMM, and SLEP on 8 selected UCI data sets with simulated groups. The regularization parameters are specified as in ($S_3$).

| problem name (m,n);g | $\gamma$ | nnz | iteration S\|A\|P | time S\|A\|P | error $\eta_S\|\eta_A\|\eta_P$ |
|---|---|---|---|---|---|
| E2006.train | 1e-05 | 1 | 4(9) \| 34 \| 16 | 01 \| 06:57 \| 01 | 8.9e-10 \| 7.8e-07 \| 8.3e-07 |
| (16087,150360) | 1e-06 | 27 | 22(45) \| 69 \| 10000 | 03 \| 07:52 \| 10:14 | 1.7e-07 \| 8.7e-07 \| 3.8e-03 |
| 5012 | 1e-07 | 1399 | 30(79) \| 395 \| 10000 | 01:27 \| 11:23 \| 09:48 | 8.5e-07 \| 9.8e-07 \| 8.6e-02 |
| E2006.test | 1e-05 | 1 | 4(10) \| 29 \| 12 | 00 \| 16 \| 00 | 8.8e-09 \| 8.8e-07 \| 7.4e-07 |
| (3308,150358) | 1e-06 | 48 | 25(51) \| 182 \| 10000 | 01 \| 25 \| 02:26 | 5.1e-08 \| 9.5e-07 \| 2.3e-02 |
| 5012 | 1e-07 | 1325 | 38(103) \| 1432 \| 10000 | 31 \| 01:17 \| 02:42 | 8.7e-09 \| 8.6e-07 \| 3.3e-01 |
| log1p.E2006.train | 1e-03 | 5 | 4(17) \| 2451 \| 626 | 13 \| 58:57 \| 04:11 | 7.8e-09 \| 9.8e-07 \| -7.5e-06 |
| (16087,4272227) | 1e-04 | 510 | 5(26) \| 823 \| 10000 | 29 \| 31:05 \| 01:01:43 | 5.5e-07 \| 9.9e-07 \| 1.1e-02 |
| 142408 | 1e-05 | 9772 | 7(33) \| 340 \| 10000 | 04:37 \| 22:48 \| 01:02:05 | 3.6e-08 \| 1.0e-06 \| 1.2e-02 |
| log1p.E2006.test | 1e-03 | 8 | 5(22) \| 1688 \| 732 | 12 \| 10:22 \| 02:09 | 5.0e-09 \| 9.9e-07 \| -3.1e-05 |
| (3308,4272226) | 1e-04 | 909 | 6(27) \| 470 \| 10000 | 18 \| 03:25 \| 29:31 | 1.0e-07 \| 9.8e-07 \| 5.4e-02 |
| 142408 | 1e-05 | 4956 | 9(35) \| 288 \| 10000 | 44 \| 02:22 \| 30:13 | 6.6e-08 \| 9.6e-07 \| 1.2e-01 |
| bodyfat7 | 1e-04 | 3 | 12(34) \| 1101 \| 1163 | 02 \| 34 \| 28 | 2.4e-07 \| 9.9e-07 \| 9.5e-07 |
| (252,116280) | 1e-05 | 24 | 19(58) \| 1586 \| 10000 | 05 \| 49 \| 05:16 | 8.6e-07 \| 1.0e-06 \| 2.1e-06 |
| 3876 | 1e-06 | 106 | 25(94) \| 2231 \| 10000 | 12 \| 01:09 \| 03:39 | 5.2e-07 \| 1.0e-06 \| 4.7e-03 |
| pyrim5 | 1e-02 | 87 | 9(32) \| 1382 \| 10000 | 01 \| 55 \| 03:51 | 4.3e-08 \| 1.0e-06 \| 7.1e-02 |
| (74,201376) | 1e-03 | 176 | 16(56) \| 5642 \| 10000 | 04 \| 03:42 \| 03:46 | 5.6e-07 \| 1.0e-06 \| 5.8e-03 |
| 6713 | 1e-04 | 129 | 26(95) \| 10000 \| 10000 | 09 \| 06:31 \| 03:52 | 5.3e-07 \| 4.1e-05 \| 1.2e-03 |
| triazines4 | 1e-02 | 246 | 10(37) \| 8369 \| 10000 | 09 \| 43:00 \| 26:58 | 5.1e-08 \| 9.2e-07 \| 6.6e-02 |
| (186,635376) | 1e-03 | 803 | 20(72) \| 10000 \| 10000 | 23 \| 50:09 \| 27:06 | 8.2e-09 \| 3.7e-06 \| 3.4e-02 |
| 21179 | 1e-04 | 333 | 27(115) \| 10000 \| 10000 | 01:04 \| 46:21 \| 27:21 | 3.9e-07 \| 1.4e-04 \| 5.5e-02 |
| housing7 | 1e-02 | 50 | 7(30) \| 976 \| 10000 | 01 \| 30 \| 04:07 | 1.3e-07 \| 1.0e-06 \| 1.2e-02 |
| (506,77520) | 1e-03 | 157 | 11(41) \| 620 \| 10000 | 03 \| 19 \| 04:03 | 1.1e-07 \| 9.9e-07 \| 6.7e-02 |
| 2584 | 1e-04 | 838 | 16(51) \| 685 \| 10000 | 08 \| 21 \| 04:11 | 8.7e-08 \| 1.0e-06 \| 3.8e-02 |

# 6 Conclusion

In this paper, we have developed a highly efficient semismooth Newton based augmented Lagrangian method SSNAL for solving large-scale non-overlapping sparse group Lasso problems. The elements in the generalized Jacobian of the proximal mapping associated with the sparse group Lasso regularizer were first derived, and the underlying second order sparsity structure was thoroughly analysed and utilised to achieve superior performance in the numerical implementations of SSNAL. Extensive numerical experiments have demonstrated that the proposed algorithm is highly efficient and robust, even on high-dimensional real data sets. Based on the superior performance of SSNAL for solving non-overlapping sparse group Lasso problems, we can expect the effectiveness of our algorithmic framework for solving overlapping sparse group Lasso problems and other large-scale convex composite problems in future studies.

# Acknowledgments

The authors would like to thank Dr. Xudong Li and Ms. Meixia Lin for their help in the numerical implementations.